\pdfoutput=1
\RequirePackage{ifpdf}
\ifpdf 
\documentclass[pdftex]{sigma}
\else
\documentclass{sigma}
\fi

\numberwithin{equation}{section}

\newtheorem{Theorem}{Theorem}[section]
\newtheorem*{Theorem*}{Theorem}

\newtheorem{Proposition}[Theorem]{Proposition}
\newtheorem{conj}[Theorem]{Conjecture}

\theoremstyle{definition}
\newtheorem{Definition}[Theorem]{Definition}

\newtheorem{Remark}[Theorem]{Remark}

\begin{document}

\allowdisplaybreaks

\newcommand{\arXivNumber}{2504.14273}

\renewcommand{\PaperNumber}{023}

\FirstPageHeading

\ShortArticleName{Elliptic Virtual Structure Constants and Gromov--Witten Invariants}

\ArticleName{Elliptic Virtual Structure Constants\\ and Gromov--Witten Invariants for Complete\\ Intersections in Weighted Projective Space}

\Author{Masao JINZENJI~$^{\rm a}$ and Ken KUWATA~$^{\rm b}$}

\AuthorNameForHeading{M.~Jinzenji and K.~Kuwata}

\Address{$^{\rm a)}$~Department of Mathematics, Okayama University, Okayama, 700-8530, Japan}
\EmailD{\mail{pcj70e4e@okayama-u.ac.jp}}

\Address{$^{\rm b)}$~Department of General Education, National Institute of Technology, Kagawa College,\\
\hphantom{$^{\rm b)}$}~Chokushi, Takamatsu, 761-8058, Japan}
\EmailD{\mail{kuwata-k@t.kagawa-nct.ac.jp}}

\ArticleDates{Received May 01, 2025, in final form February 17, 2026; Published online March 09, 2026}

\Abstract{In this paper, we generalize our formalism of the elliptic virtual structure constants to hypersurfaces and complete intersections within certain weighted projective spaces possessing a single K\"{a}hler class.}


\Classification{14N35; 81T45}

\section{Introduction}
In this paper, we generalize the definition of elliptic virtual structure constants for projective hypersurfaces, as presented in \cite{MK}, to encompass hypersurfaces and complete intersections in specific weighted projective spaces with one K\"{a}hler class. We then compute their genus~1 (elliptic) Gromov--Witten invariants. The guideline for this generalization stems from previous works~\cite{M2,MS} of our group, which detail the construction of genus 0 virtual structure constants for the K3 surface in the weighted projective space $P(1,1,1,3)$.
However, akin to the situation in~\cite{MK}, our definition lacks a rigorous geometrical construction of the expected moduli space of quasimaps from an elliptic curve to weighted projective spaces. Consequently, we only explicitly write down the integrands of the residue integrals associated with the graph types~(i),~(ii),~(iii), and~(iv) introduced in \cite{MK} (the necessary graph types remain the same as for projective hypersurfaces). A primary motivation for this generalization is to validate our conjecture by comparing the results obtained from our formalism for certain complex $3$-dimensional Calabi--Yau hypersurfaces in weighted projective spaces with one K\"{a}hler class against the corresponding results derived from the original BCOV formalism \cite{BCOV}.
Let $P(a_1,a_2,\dots,a_N\mid k_1,k_2,\dots,k_m)$ denote a~complete intersection of degree $k_{1},\dots,k_{m}$ homogeneous polynomials within the weighted projective space
$P(a_{1},a_{2},\dots,a_{N})$ with a single K\"{a}hler class. The significant findings resulting from this generalization can be summarized in the following two points:
\begin{itemize}\itemsep=0pt
\item[(i)] The nontrivial part:
\[
\left(-\frac{N-1}{N}\frac{1}{w^N}-\frac{N+1}{N}\frac{1}{(z_0)^N} \right),
\]
 which appeared in the integrand associated with the type (iii) graph (see \cite{MK}), is
 modified to
 \[
 \left(-\frac{N-m}{N}\frac{1}{w^N}-\frac{N+m}{N}\frac{1}{(z_0)^N} \right).
 \]
 \item[(ii)] The symmetric factor:
\[
R_{N,k}(d):=\left(\frac{N-1}{2}\right)\frac{1}{d}-\left(N-\frac{1}{k}\right)\frac{1}{d^2},
\]
 in \cite{MK} associated with the type (iv) graph is modified to
\[
R(d):=\left(\prod_{i=1}^Na_i\right)\left( \left(\frac{N-m}{2}\right) \frac{1}{d}-\left( \sum_{j=1}^N \frac{1}{a_j}-\sum_{l=1}^m \frac{1}{k_l}\right)\frac{1}{d^2}\right).
\]
\end{itemize}
The remaining generalizations follow straightforwardly from the results presented in \cite{MK} and our prior findings concerning the K3 surface in the weighted projective space $P(1,1,1,3)$ \cite{M2,MS}.

This paper unfolds as follows:
First, in Section~\ref{Complete Intersection in Weighted Projective Space}, we lay the groundwork by reviewing essential properties of weighted projective spaces and the complete intersections they contain, drawing upon the established literature \cite{AR, AS, PV}.
Moving on to Section~\ref{Set up of Our Computation}, we then elucidate the theoretical framework we employ -- our formalism detailed in \cite{M2} -- for the computation of Gromov--Witten invariants through the lens of virtual structure constants.
Subsequently, Section~\ref{Elliptic Virtual Structure Constants} provides a precise definition of the elliptic virtual structure constant. This definition takes the form of a sum of residue integrals, where each integrand corresponds to one of the four graph types originally introduced in~\cite{MK}.
In Section~\ref{Numerical Tests for Various Examples}, we put our formalism to the test by presenting a range of numerical results. We begin by examining the genus~1 Gromov--Witten invariants of Fano hypersurfaces within specific weighted projective spaces. Following this, we~turn our attention to Calabi--Yau 3-folds residing in the weighted projective space previously studied in~\cite{BCOV}. Notably, our findings regarding the count of elliptic (and rational) curves for these Calabi--Yau manifolds are consistent with those reported in~\cite{BCOV}. We then extend our analysis to complete intersections within standard projective space and, finally, to complete intersections in various weighted projective spaces.
Lastly, for the interested reader, Appendix~\ref{appendixA} offers a formal proof of a proposition concerning the vanishing of residue integrals associated with type (iii) graphs, specifically in the context of Calabi--Yau manifolds.

\section{Complete intersection in weighted projective space}\label{Complete Intersection in Weighted Projective Space}

In this section, we summarize the properties of weighted projective space and complete intersections within it. First, we introduce the weighted projective space $P(a_1, a_2, \dots, a_N)$, where $a_1, a_2, \dots, a_N$ are positive integers.

From this point onward, we restrict $a_1$ to be $1$ and assume that the greatest common divisor of $a_2, a_3, \dots, a_N$ is $1$. Let $\bigl(x^1, x^2, \dots, x^N\bigr)$ be the coordinates on $\mathbb{C}^N$. We define an action of the multiplicative group of complex numbers, $\mathbb{C}^{\times}$, on $\mathbb{C}^N \setminus \{\mathbf{0}\}$ as follows:
\begin{align*}
\lambda \cdot \bigl(x^1, x^2, \dots, x^N\bigr) = \bigl(\lambda x^1, \lambda^{a_2} x^2, \dots, \lambda^{a_N} x^N\bigr), \qquad \lambda \in \mathbb{C}^{\times}.
\end{align*}
Then, ${P}(a_1, a_2, \dots, a_N)$ is defined as the orbit space of this $\mathbb{C}^{\times}$ action,
\begin{align*}
{P}(a_1, a_2, \dots, a_N) = \bigl(\mathbb{C}^N \setminus \{\mathbf{0}\}\bigr) / \mathbb{C}^{\times}.
\end{align*}
We denote a point in $P(a_1, a_2, \dots, a_N)$ represented by $(x_1, x_2, \dots, x_N) \in \mathbb{C}^N \setminus \{\mathbf{0}\}$ as $(x_1 : x_2 : \dots : x_N)$.

In general, $P(a_1, a_2, \dots, a_N)$ can have singularities. Information about its singular locus can be obtained from the toric construction of the weighted projective space. It is constructed as a~toric variety associated with a polytope $\Delta^* \subset \mathbb{R}^{N-1}$. The polytope~$\Delta^*$ is given as the convex hull of the following $N$ points in $\mathbb{R}^{N-1}$:
\begin{align*}
v_1 = (-a_2, -a_3, \dots, -a_N), \qquad v_i = e_{i-1}, \quad i = 2, \dots, N,
\end{align*}
where $e_{i-1}$ is the $(i-1)$-th standard basis vector of $\mathbb{R}^{N-1}$. Let $\langle v_{i_1}, \dots, v_{i_m} \rangle_{\mathbb{R}_{\geq 0}}$, $1 \leq i_1 < i_2 < \dots < i_m \leq N$, be the $m$-dimensional cone in $\mathbb{R}^{N-1}$ defined as follows:
\begin{align*}
\langle v_{i_1}, \dots, v_{i_m} \rangle_{\mathbb{R}_{\geq 0}} := \Biggl\{ \sum_{j=1}^{m} u_j v_{i_j} \in \mathbb{R}^{N-1} \mid u_1, \dots, u_m \geq 0 \Biggr\}.
\end{align*}
The conventional rule for identifying the singular locus is as follows:
\begin{itemize}\itemsep=0pt
 \item[(i)] The locus $\{\;(x_1:\cdots:x_N)\in P(a_{1},\dots,a_N)\;|\; x_{i_{1}}=x_{i_{2}}=\cdots=x_{i_{m}}=0\}$ is singular if and only if the cone $\langle v_{i_{1}},\cdots v_{i_{m}} \rangle_{\mathbb{R}_{\geq 0}}$ contains integral points that cannot be expressed as an integral linear combination of $v_{i_{1}},v_{i_{2}},\dots,v_{i_{m}}$.
\end{itemize}

Let us consider ${P}(1,1,1,2)$ as an example. The corresponding polytope $\Delta^{*}$ is the convex hull of the following four points in $\mathbb{R}^{3}$:
\begin{align*}
v_{1}=(-1,-1,-2), \qquad v_{2}=(1,0,0), \qquad v_{3}=(0,1,0), \qquad v_{4}=(0,0,1).
\end{align*}
We can easily observe that only the cone $\langle v_{1},v_{2},v_{3} \rangle_{\mathbb{R}_{\geq 0}}$ contains integral points that cannot be expressed as an integral linear combination of $v_{1},v_{2}$, and $v_{3}$ (for example, $(0,0,-1)$). Therefore, the singular locus of ${P}(1,1,1,2)$ is given by the point $(0:0:0:1) \in {P}(1,1,1,2)$. The singular locus of the weighted projective spaces treated in this paper is always a single point.

On the other hand, we have the following exact sequence:
\begin{align*}
\mathbf{0} \to \mathbb{C}^{N-1} \to \bigoplus_{j=1}^{N} \mathbb{C} v_j \to CH^{1}({P}(a_{1},a_{2},\dots,a_{N})) \to \mathbf{0},
\end{align*}
where $CH^{1}(M)$ is the degree $1$ Chow ring of an algebraic variety $M$ (if $M$ is non-singular, it~equals $H^{1,1}(M,\mathbb{C})$). Hence, $\dim_{\mathbb{C}}\bigl(CH^{1}({P}(a_{1},a_{2},\dots,a_{N}))\bigr)$ is $1$. In the cases of hypersurfaces and complete intersections in the weighted projective spaces treated in this paper, we can make them non-singular by avoiding the singular points of the weighted projective spaces through an appropriate choice of defining equations. Therefore, the resulting hypersurface or complete intersection~$M$ becomes non-singular and has a single K\"{a}hler class.

Next, we introduce the concept of a non-singular complete intersection in the weighted projective space ${P}(a_1, a_2, \dots, a_N)$. Let $F_1, F_2, \dots, F_m$ be weighted homogeneous polynomials with degrees $k_1, k_2, \dots, k_m$, respectively. The degree $k_i$ hypersurface in the weighted projective space is defined as the zero locus of $F_i$ in ${P}(a_1, a_2, \dots, a_N)$, and we denote it by ${P}(a_1, a_2, \dots, a_N\mid k_i)$.

In what follows, we assume that the hypersurface is well-formed and quasi-smooth. In particular, when it is Cartier, its degree is divisible by the relevant weights. If some $a_i$ coincides with~$k_1$, then for a general weighted complete intersection the corresponding factor can be eliminated without affecting
the Gromov--Witten invariants.

The complete intersection $P(a_1, a_2, \dots, a_N\mid k_1, k_2, \dots, k_m)$ is defined by
\begin{align*}
P(a_1, a_2, \dots, a_N\mid k_1, k_2, \dots, k_m) = P(a_1, a_2, \dots, a_N\mid k_1) \cap \cdots \cap P(a_1, a_2, \dots, a_N\mid k_m)
\end{align*}
where $k_i$, $i = 1, \dots, m$, are positive integers and $a_j$ divides~$k_i$ for all $i$ and~$j$.

Lastly, we introduce the Calabi--Yau threefolds treated in this paper. In Section~\ref{Numerical Tests for Various Examples}, we compute the genus $1$ Gromov--Witten invariants for the following three types of Calabi--Yau hypersurfaces, which were also used as examples in~\cite{AS}:
\begin{align*}
&P(1,1,1,1,2\mid 6) = \left \{ \sum_{i=1}^4 \bigl(x^i\bigr)^6 + 2\bigl(x^5\bigr)^3 = 0 \right \}, \\
&P(1,1,1,1,4\mid 8) = \left \{ \sum_{i=1}^4 \bigl(x^i\bigr)^8 + 4\bigl(x^5\bigr)^2 = 0 \right \}, \\
&P(1,1,1,2,5\mid 10) = \left \{ \sum_{i=1}^3 \bigl(x^i\bigr)^{10} + 2\bigl(x^4\bigr)^5 + 5\bigl(x^5\bigr)^2 = 0 \right \}.
\end{align*}

\section{Setup of our computation}\label{Set up of Our Computation}

In this section, we outline our strategy to compute genus $1$ Gromov--Witten invariants of the non-singular complete intersection $P(a_1,\dots,a_N\mid k_1,\dots,k_m)$ in the weighted projective space $P(a_1,\dots,a_N)$.
Let $h$ be the hyperplane (K\"{a}hler) class of $P(a_1,\dots,a_N)$. We also denote the restriction of this hyperplane class to the complete intersection by the same symbol $h$.
We~then introduce the genus $g$ Gromov--Witten invariants
\[
\biggl\langle\prod_{a=2}^{N-m-1} ( {\mathcal O}_{h^a})^{n_a} \biggr\rangle_{g,d}
\]
 of~$P(a_1,\dots,a_N\mid k_1,\dots,k_m)$, which represent the intersection number of the moduli space of holomorphic maps from genus~$g$ stable curves $\Sigma_{g}$ with $\sum_{a=2}^{N-m-1}n_{a}$ marked points to $P(a_1,a_2,\allowbreak\dots, a_N\mid k_1,k_2,\dots,k_m)$ of degree~$d$. We denote this moduli space by
\[
\overline{M}_{g,\sum_{a=2}^{N-m-1}n_{a}}(P(a_1,a_2,\dots,a_N\mid k_1,k_2,\dots,k_m),d).
\]
Here, ${\mathcal O}_{h^{a}}$ represents the insertion of the pullback of the cohomology class $h^{a}$ by the evaluation map
\begin{gather*}
\mathrm{ev}_{i}\colon \ \overline{M}_{g,\sum_{a=2}^{N-m-1}n_{a}}(P(a_1,a_2,\dots,a_N\mid k_1,k_2,\dots,k_m),d)\\
\hphantom{\mathrm{ev}_{i}\colon} \ \rightarrow P(a_1,a_2,\dots,a_N\mid k_1,k_2,\dots,k_m),
\end{gather*}
at the $i$-th marked point $z_{i}\in\Sigma_{g}$, $i=1,2,\dots,\sum_{a=2}^{N-m-1}n_{a}$. The complex dimension of the moduli space is given by
\[
N-m-1+ \left(\sum_{i=1}^Na_i-\sum_{j=1}^mk_j\right)d+(3g-3)+\sum_{a=2}^{N-m-1}n_{a}.
\]
 Therefore, $\bigl\langle\prod_{a=2}^{N-m-1} ( {\mathcal O}_{h^a})^{n_a} \bigr\rangle_{g,d}$ is non-zero only if the
following condition is satisfied:
\begin{align}
&N-m-1+ \left(\sum_{i=1}^Na_i-\sum_{j=1}^mk_j\right)d+(3g-3)+\sum_{a=2}^{N-m-1}n_{a}=\sum_{a=2}^{N-m-1}n_{a}a \nonumber\\
&\qquad \Longleftrightarrow \ N-m-1+ \left(\sum_{i=1}^Na_i-\sum_{j=1}^mk_j\right)d+(3g-3)=\sum_{a=2}^{N-m-1}n_{a}(a-1).
\label{selgw}
\end{align}
Next, we introduce the genus $0$ multi-point virtual structure constants
\[
w\Biggl({\mathcal O}_{h^a} {\mathcal O}_{h^b} \,\Big|\, \prod_{j=0}^{N-m-1} ( {\mathcal O}_{h^j})^{n_j} \Biggr)_{0,d}
\]
 of the same complete intersection. To simplify the notation, we define the following polynomials:
\begin{align*}
&e_k(x,y):=\prod_{i=0}^k(ix+(k-i)y),\qquad 
w_a(x,y):=\frac{x^a-y^a}{x-y}, \\ 
&q(x,y):=\prod_{i=1}^N\prod_{j=1}^{a_i-1}(jx+(a_i-j)y).
\end{align*}
If $a_i=1$ for some $i$, then the term $\prod_{j=1}^{a_i-1}(jx+(a_i-j)y)$ is defined to be $1$. With these definitions, the multi-point virtual structure constant for $d\geq 1$ is given by
\begin{align*}
&w\Biggl({\mathcal O}_{h^a} {\mathcal O}_{h^b}\,\Big|\, \prod_{p=0}^{N-m-1} ( {\mathcal O}_{h^p})^{n_p} \Biggr)_{0,d} \notag \\
&\qquad =\frac{1}{(2\pi \sqrt{-1})^{d+1}} \oint _{C_{z_0}}\mathrm{d}z_0 \oint _{C_{z_1}}\mathrm{d} z_1 \cdots \oint _{C_{z_d}}\mathrm{d}z_d \frac{(z_0)^a (z_d)^b}{\bigl(\prod_{i=1}^N a_i\bigr)^{d+1}\prod_{i=0}^d (z_i)^N}\notag \\
&\qquad\quad{} \times \prod_{j=1}^{d} \left( \frac{\prod_{p=1}^{m} e_{k_p}(z_{j-1}, z_j)}{\bigl( \prod_{p=1}^{m}k_p z_j\bigr)(2z_j-z_{j-1}-z_{j+1}) }\right) \Biggl( \prod_{p=1}^{m}k_p z_d\Biggr)
\frac{1}{\prod_{j=1}^d q(z_{j-1}, z_j)}\notag \\
&\qquad\quad{} \times\left( \prod_{l=0}^{N-m-1}\left(\sum_{i=1}^d w_l(z_{i-1},z_{i})
 \right) ^{n_l}\right).
\end{align*}
Here, $\frac{1}{2\pi\sqrt{-1}}\oint_{C_{z_{i}}}\mathrm{d}z_i$ represents taking the residue at $z_i=0$ if $i=0$ or $i=d$, and at $z_i=\frac{z_{i-1}+z_{i+1}}{2}$ or $z_i=0$ if $i=1,2,\dots,d-1$.
From this definition, we can see that it obeys the following selection rule:
\begin{align*}
&w\Biggl({\mathcal O}_{h^a} {\mathcal O}_{h^b}\,\Big|\, \prod_{j=0}^{N-m-1} ( {\mathcal O}_{h^j})^{n_j} \Biggr)_{0,d} \neq 0 \notag \\
&\qquad \Rightarrow \ \left(\sum_{i=1}^Na_i-\sum_{j=1}^mk_j\right)d+N-m-2=a+b+\sum_{j=0}^{N-m-1}n_j(j-1).
\end{align*}
Since $w_{0}(x,y)=0$ and $w_1(x,y)=1$, the following condition follows:
\begin{align}
w\Biggl({\mathcal O}_{h^a} {\mathcal O}_{h^b}\,\Big|\, \prod_{p=0}^{N-m-1} ( {\mathcal O}_{h^p})^{n_p} \Biggr)_{0,d} =\delta_{n_{0},0}\cdot d^{n_1} w\Biggl({\mathcal O}_{h^a} {\mathcal O}_{h^b}\,\Big|\, \prod_{p=2}^{N-m-1} ( {\mathcal O}_{h^p})^{n_p} \Biggr)_{0,d}.
\label{wch1}
\end{align}
If $d=0$, the multi-point virtual constants vanish except for the following cases:
\begin{align}
w({\mathcal O}_{h^a} {\mathcal O}_{h^b}\mid {\mathcal O}_{h^c} )_{0,0} = \frac{\prod_{l=1}^mk_l}{\prod_{i=1}^Na_i} \delta_{a+b+c, N-m-1}.\label{wch2}
\end{align}
We introduce two types of perturbed two-point functions: $w({\mathcal O}_{h^a} {\mathcal O}_{h^b})_{0}(x^*) $ and $\langle{\mathcal O}_{h^a} {\mathcal O}_{h^b}\rangle_{0}(t^*) $. Here, $t^*$ and $x^*$ denote deformation variables $t^0,t^1,t^2, \dots,t^{N-m-1}$ and $x^0,x^1,x^2,\dots,x^{N-m-1}$, respectively. These functions are generating functions of the multi-point virtual structure constants and the Gromov--Witten invariants, respectively. The first one, $w({\mathcal O}_{h^a} {\mathcal O}_{h^b})_{0}(x^*) $, is defined as follows:
\begin{align*}
&w({\mathcal O}_{h^a}{\mathcal O}_{h^b})_{0}\bigl(x^0,x^1,x^2,\dots,x^{N-m-1}\bigr) \\
&\qquad =\sum_{n_0=0}^\infty \sum_{n_1=0}^\infty \sum_{n_2=0}^\infty \cdots \sum_{n_{N-m-1}=0}^\infty \sum_{d=0}^\infty w\Biggl({\mathcal O}_{h^a} {\mathcal O}_{h^b}\,\Big|\, \prod_{q=0}^{N-m-1} ( {\mathcal O}_{h^q})^{n_q} \Biggr)_{0,d} \prod_{q=0}^{N-m-1} \frac{(x^q)^{n_q}}{n_q!}.
\end{align*}
By using the properties~(\ref{wch1}) and~(\ref{wch2}), we can rewrite the right-hand side as follows:
\begin{align}
&w({\mathcal O}_{h^a} {\mathcal O}_{h^b})_{0}\bigl(x^0,x^1,x^2,\dots,x^{N-m-1}\bigr)
 =\frac{\prod_{l=1}^mk_l}{\prod_{i=1}^Na_i}\bigl(x^1\bigr)^{N-m-1-a-b}\nonumber\\
&\qquad{} +\sum_{n_2=0}^\infty \sum_{n_3=0}^\infty \cdots \sum_{n_{N-m-1}=0}^\infty \sum_{d=1}^\infty e^{dx^1}w\Biggl({\mathcal O}_{h^a} {\mathcal O}_{h^b}\,\Big|\, \prod_{q=2}^{N-m-1} ( {\mathcal O}_{h^q})^{n_q} \Biggr)_{0,d} \prod_{q=2}^{N-m-1} \frac{(x^q)^{n_q}}{n_q!}.\label{wstr}
\end{align}
The second one, $\langle{\mathcal O}_{h^a} {\mathcal O}_{h^b}\rangle_{0}(t^*) $, is defined by
\begin{align*}
&\langle{\mathcal O}_{h^a} {\mathcal O}_{h^b} \rangle_{0}\bigl(t^0,t^1,t^2,\dots,t^{N-m-1}\bigr)
 =\frac{\prod_{l=1}^mk_l}{\prod_{i=1}^Na_i}\bigl(t^1\bigr)^{N-m-1-a-b} \\
 &\qquad {} +\sum_{n_2=0}^\infty \sum_{n_3=0}^\infty \cdots \sum_{n_{N-m-1}=0}^\infty \sum_{d=1}^\infty e^{dt^1} \Biggl\langle{\mathcal O}_{h^a} {\mathcal O}_{h^b} \prod_{q=2}^{N-m-1} ( {\mathcal O}_{h^q})^{n_q} \Biggr\rangle_{0,d} \prod_{q=2}^{N-m-1} \frac{(t^q)^{n_q}}{n_q!}.
\end{align*}
At this stage, we rewrite (\ref{selgw}), the topological selection rule for the Gromov--Witten invariants of $P(a_1,\dots,a_N\mid k_1,\dots,k_m)$ for the specific $g=0,1$ cases:
\begin{align*}
&\Biggl\langle \prod_{a=0}^{N-m-1} ( {\mathcal O}_{h^a})^{n_a} \Biggr\rangle_{0,d} \neq 0 \ \Rightarrow \ \Biggl(\sum_{i=1}^Na_i-\sum_{j=1}^mk_j\Biggr)d+N-m-2=\sum_{j=0}^{N-m-1}n_j(j-1),\\
&\Biggl\langle\prod_{a=0}^{N-m-1} ( {\mathcal O}_{h^a})^{n_a} \Biggr\rangle_{1,d} \neq 0 \ \Rightarrow \ \Biggl(\sum_{i=1}^Na_i-\sum_{j=1}^mk_j \Biggr)d=\sum_{j=0}^{N-m-1}n_j(j-1).
\end{align*}
With these setups, we define the mirror maps $t^p(x^*)$, $p=0,1,2,\dots,N-m-1$, using the first type of the two-point perturbed function as follows:
\begin{align*}
&t^p\bigl(x^0,x^1,\dots,x^{N-m-1}\bigr):=\frac{\prod_{i=1}^Na_i}{\prod_{l=1}^{m}k_l}w({\mathcal O}_{h^{N-m-1-p}} {\mathcal O}_{1})_{0} (x^*), \qquad p=0,1,2,\dots,N-m-1 .
\end{align*}
We can see from (\ref{wstr}) that it has the following structure:
\begin{align*}
t^p\bigl(x^0,x^1,\dots,x^{N-m-1}\bigr)
 ={}&x^{p}+\frac{\prod_{i=1}^Na_i}{\prod_{l=1}^{m}k_l}\sum_{n_2=0}^\infty \sum_{n_3=0}^\infty \cdots \sum_{n_{N-m-1}=0}^\infty \sum_{d=1}^\infty e^{dx^1}\\
 & \times w\Biggl({\mathcal O}_{h^{N-m-1-p}} {\mathcal O}_{1}\,\Big|\, \prod_{q=2}^{N-m-1} ( {\mathcal O}_{h^q})^{n_q} \Biggr)_{0,d} \prod_{q=2}^{N-m-1} \frac{(x^q)^{n_q}}{n_q!},
\end{align*}
which enables us to invert the mirror map
\begin{align*}
 x^p=x^p\bigl(t^0,t^1,\dots,t^{N-m-1}\bigr), \qquad p=0,1,2,\dots,N-m-1.
\end{align*}
In this paper, we use the following conjecture~\cite{M3} to compute the genus~$0$ Gromov--Witten invariants, which has already been proved in~\cite{M6} for the case of projective hypersurfaces.

\begin{conj}
\begin{align*}
&\langle{\mathcal O}_{h^a} {\mathcal O}_{h^b}\rangle_{0}\bigl(t^0,t^1,t^2,\dots,t^{N-m-1}\bigr)=w({\mathcal O}_{h^a} {\mathcal O}_{h^b})_{0}\bigl(x^0(t^*),x^1(t^*),x^2(t^*),\dots,x^{N-m-1}(t^*)\bigr).
\end{align*}
\end{conj}

Lastly, we introduce the generating function for genus $1$ Gromov--Witten invariants. In this paper, we also introduce the generating function of the elliptic virtual structure constants $w\bigl( \prod_{p=0}^{N-m-1} ( {\mathcal O}_{h^p})^{n_p} \bigr)_{1,d}$, whose general definition will be given in the next section, to compute the genus $1$ Gromov--Witten invariants.
In the $d=0$ case, $w\bigl( \prod_{p=0}^{N-m-1} ( {\mathcal O}_{h^p})^{n_p} \bigr)_{1,0}$ vanishes except for the following case:
\begin{align*}
w( {\mathcal O}_{h} )_{1,0}:=\langle {\mathcal O}_{h} \rangle_{1,0}& =-\frac{1}{24}\int_M c_{N-m-2}(P(a_1,\dots,a_N\mid k_1,\dots,k_m))h\\
& =-\frac{\prod_{l=1}^{m}k_l}{24\prod_{j=1}^N a_j}c_{N-m-2}.
\end{align*}
Here, $h$ is the K\"{a}hler class of $P(a_1,\dots,a_N\mid k_1,\dots,k_m)$. $c_{N-m-2}( P(a_1,\dots,a_N\mid k_1,\dots,k_m) )$ is the second top Chern class of $P(a_1,a_2,\dots,a_N\mid k_1,k_2,\dots,k_m)$. $c_{N-m-2}$ in the right-hand side is given by the following formula:
\begin{equation*}
c(P(a_1,\dots,a_N\mid k_1,\dots,k_m))=\frac{\prod_{j=1}^N(1+ta_jh)}{\prod_{l=1}^m(1+tk_lh)}=:\sum_{j=0}^{N-m-1}c_{j}h^{j}.
\end{equation*}
If $d\geq 1$, $w\bigl( \prod_{p=0}^{N-m-1} ( {\mathcal O}_{h^p})^{n_p} \bigr)_{1,d}$ obeys the following conditions:
\begin{itemize}\itemsep=0pt
 \item[(i)] $w\bigl(\prod_{p=0}^{N-m-1} ( {\mathcal O}_{h^p})^{n_p}\bigr)_{1,d} \neq 0 \Rightarrow \bigl(\sum_{i=1}^Na_i-\sum_{j=1}^mk_j \bigr)d=\sum_{j=0}^{N-m-1}n_j(j-1).$
 \item[(ii)] $w\bigl( \prod_{p=0}^{N-m-1} ( {\mathcal O}_{h^p})^{n_p} \bigr)_{1,d}=\delta_{n_0,0} d^{n_1} w\bigl( \prod_{p=2}^{N-m-1} ( {\mathcal O}_{h^p})^{n_p} \bigr)_{1,d}$.
\end{itemize}
We then introduce the generating function of $w\bigl( \prod_{a=0}^{N-m-1} ( {\mathcal O}_{h^a})^{n_a} \bigr)_{1,d}$,
\begin{align*}
&F^{B}_{1}\bigl(x^0,x^1,\dots,x^{N-m-1}\bigr)\notag \\
&\qquad =\sum_{n_0=0}^\infty\sum_{n_1=0}^\infty \sum_{n_2=0}^\infty \cdots \sum_{n_{N-m-1}=0}^\infty \sum_{d=1}^\infty w\Biggl( \prod_{p=0}^{N-m-1} ( {\mathcal O}_{h^p})^{n_p} \Biggr)_{1,d} \prod_{q=0}^{N-m-1} \frac{(x^q)^{n_q}}{n_q!} \\
&\qquad =-\frac{\prod_{l=1}^{m}k_l}{24\prod_{i=1}^Na_i}c_{N-m-2}x^1\notag \\
&\qquad\quad{} +\sum_{n_2=0}^\infty \sum_{n_3=0}^\infty \cdots \sum_{n_{N-m-1}=0}^\infty \sum_{d=1}^\infty e^{dx^1}w\Biggl( \prod_{p=2}^{N-m-1} ( {\mathcal O}_{h^p} )^{n_p} \Biggr)_{1,d} \prod_{q=2}^{N-m-1} \frac{(x^q)^{n_q}}{n_q!},
\end{align*}
and the generating function of genus one Gromov--Witten invariants
\begin{align*}
&F^{A}_{1}\bigl(t^0,t^1,\dots,t^{N-m-1}\bigr)\notag \\
&\qquad =\sum_{n_0=0}^\infty \sum_{n_1=0}^\infty \cdots \sum_{n_{N-m-1}=0}^\infty \sum_{d=1}^\infty \Biggl\langle\prod_{p=0}^{N-m-1} ( {\mathcal O}_{h^p})^{n_p} \Biggr\rangle_{1,d} \prod_{q=0}^{N-m-1} \frac{(t^q)^{n_q}}{n_q!}. \\
&\qquad =-\frac{\prod_{l=1}^{m}k_l}{24\prod_{i=1}^Na_i}c_{N-m-2}t^1\notag \\
&\qquad \quad{} +\sum_{n_2=0}^\infty \sum_{n_3=0}^\infty \cdots \sum_{n_{N-m-1}=0}^\infty \sum_{d=1}^\infty e^{dt^1}\Biggl\langle\prod_{p=2}^{N-m-1} ( {\mathcal O}_{h^p})^{n_p} \Biggr\rangle_{1,d} \prod_{q=2}^{N-m-1} \frac{(t^q)^{n_q}}{n_q!}.
\end{align*}
Our conjecture used for computing genus $1$ Gromov--Witten invariants of the non-singular complete intersection $P(a_1,\dots,a_N \mid k_1,\dots,k_m)$ is given as follows.

\begin{conj}\label{main}
\begin{align*}
&F^{A}_{1}\bigl(t^0,t^1,\dots,t^{N-m-1}\bigr)=F^{B}_{1}\bigl(x^0(t^*),x^1(t^*),x^2(t^*),\dots,x^{N-m-1}(t^*)\bigr).
\end{align*}
\end{conj}

\begin{Remark}
In the case when $a_1=a_2=\cdots=a_N=1$ and $m=1$, these expressions reduce to the ones used in our previous paper \cite{MK}.
\end{Remark}

\section{Elliptic virtual structure constants} \label{Elliptic Virtual Structure Constants}

\subsection[Graphs for elliptic virtual structure constants \protect{[7]}]{Graphs for elliptic virtual structure constants \cite{MK}} \label{Graphs for Elliptic Virtual Structure Constants}

The elliptic virtual structure constants $w\bigl( \prod_{p=0}^{N-m-1} ( {\mathcal O}_{h^p})^{n_p} \bigr)_{1,d}$ are defined by the summation of residues of integrands associated with graphs. In this section, we briefly introduce the four types of graphs used in their definition. For more detailed explanations or figures of the graphs, we recommend that readers refer to~\cite{MK}.

For preparation, we introduce the partition of a positive integer~$d$,
\[
 \sigma=(d_1,d_2,\dots,d_{l(\sigma)}), \qquad d_1\leq d_2\leq d_3\leq \cdots \leq d_{l(\sigma)}, \qquad \sum_{i=1}^{l(\sigma)} d_i=d.
 \]
We call $l(\sigma)$ the length of the partition $\sigma$. Let $P_d$ be the set of partitions of the positive integer~$d$:%
\[
 P_d:=\Biggl\{\sigma=(d_1,d_2,\dots,d_{l(\sigma)}) \mid d_1\leq d_2\leq d_3\leq \cdots \leq d_{l(\sigma)},\, \sum_{i=1}^{l(\sigma)} d_i =d\Biggr\}.
 \]
We then define the symmetry factor associated with $\sigma \in P_d$:
\[
\mathrm{Sym}(\sigma)=\frac{(l(\sigma)-1)!}{\prod_{i=1}^{l(\sigma)} \mathrm{mul}(\sigma,i)!}, \qquad \sigma \in P_d.
\]
Here, $\mathrm{mul}(\sigma,i)$ denotes the multiplicity of $i$, $1\leq i \leq d$, in $\sigma$.

The graphs used in the definition of $w\bigl( \prod_{p=0}^{N-m-1} ( {\mathcal O}_{h^p})^{n_p} \bigr)_{1,d}$ are constructed using edges and three types of vertices:
\begin{itemize}\itemsep=0pt
 \item[(i)] normal vertex,
 \item[(ii)] elliptic vertex,
 \item[(iii)] cluster vertex of degree $d$.
\end{itemize}
A single edge is assigned a degree of $1$. A normal vertex, an elliptic vertex, and a cluster vertex of degree $d$ are assigned degrees of $0$, $0$, and $d$, respectively.

Next, we introduce the following four types of graphs:
\begin{itemize}\itemsep=0pt
 \item[(i)] star graph associated with $\sigma \in P_d$,
 \item[(ii)] loop graph with $d$ edges and $d$ normal vertices ($d\ge 2$),
 \item[(iii)] star graph associated with $\sigma \in P_{d-f}$, $1\leq f \leq d-1$, having a cluster vertex of degree~$d$ as its center,
 \item[(iv)] graph consisting of a single cluster vertex of degree~$d$.
\end{itemize}

In this paper, we denote the sets of graphs of type~(i) and~(iii) by $\mathrm{Graph}^{\rm (i)}_d$ and $\mathrm{Graph}^{\rm (iii)}_d$, respectively. We can easily see that the type~(ii) graph and the type (iv) graph of degree $d$ are unique, and we denote these graphs by $\Gamma^{\mathrm{loop}}_d$ and $\Gamma^{\mathrm{point}}_d$, respectively. We also denote by $f_\Gamma$ the integrand associated with the graph $\Gamma$, which will be defined in the next subsection. With these setups, the elliptic virtual structure constant $w\bigl( \prod_{a=0}^{N-m-1} ( {\mathcal O}_{h^a})^{n_a} \bigr)_{1,d}$ is defined as follows.

\begin{Definition}
\begin{gather*}
 w\Biggl( \prod_{a=0}^{N-m-1} ( {\mathcal O}_{h^a})^{n_a}\Biggr)_{1,d}:=\sum_{\Gamma \in \mathrm{Graph}^{\rm (i)}_d } \mathrm{Res}(f_\Gamma) +\mathrm{Res}(f_{\Gamma^{\mathrm{loop}}_d})\\
 \hphantom{w\Biggl( \prod_{a=0}^{N-m-1} ( {\mathcal O}_{h^a})^{n_a}\Biggr)_{1,d}:=}{}
 +\sum_{\Gamma \in \mathrm{Graph}^{\rm (iii)}_d } \mathrm{Res}(f_\Gamma)+\mathrm{Res}(f_{\Gamma^{\mathrm{point}}_d}).
\end{gather*}
Here, $\mathrm{Res}(f_\Gamma)$ denotes the procedure of taking the residue of $f_\Gamma$, which will also be explained in the next subsection.
\end{Definition}

\subsection{Residue integrals for elliptic virtual structure constants}

In this subsection, we explicitly write down the integrand $f_{\Gamma}$ associated with the graph $\Gamma$.
First, we recall the total Chern class of the complete intersection $P(a_1,\dots,a_N\mid k_1,\dots,k_m)$, as given in \cite{BR}:
\begin{align*}
c(P(a_1,\dots,a_N\mid k_1,\dots,k_m))=\frac{\prod_{j=1}^N(1+ta_jh)}{\prod_{l=1}^m(1+tk_lh)}=:\sum_{j=0}^{N-m-1}c_{j}h^{j}.
\end{align*}
Then, we define $c_T(z)$ as follows:
\begin{align*}
c_T(z):=c_{N-m-1}z.
\end{align*}
The integrand $f_\Gamma$ associated with a graph $\Gamma$ is defined as follows:
For a type (i) star graph of degree $d$ associated with a partition $\sigma=(d_1,d_2,\dots,d_l) \in P_{d}$, we prepare $d+1$ variables~$z_0$ and~$z_{i,j}$ (where $1\leq i \leq l$, $1\leq j \leq d_j$). Here, $z_0$~is associated with the elliptic vertex, and $z_{i, j}$ is associated with a normal vertex. Then, $f_\Gamma$ for this graph is given by
\begin{align*}
&\frac{\text{Sym}(\sigma)}{24\bigl(\prod_{i=1}^Na_i\bigr)^{d+1}(z_0)^N} \Biggl( \prod_{i=1}^l\Biggl(\prod_{j=1}^{d_i}\frac{1}{(z_{i,j})^N}
 \Biggr) \Biggr) \frac{c_T(z_0)\prod_{a=2}^{N-m-1}\bigl(\sum_{i=1}^l \sum_{j=1}^{d_i}w_a(z_{i,j-1},z_{i,j})
 \bigr) ^{n_a}}{\prod_{l=1}^m(k_lz_0)^{l-1}} \\
 &\times\Biggl( \prod_{i=1}^l\frac{\prod_{l=1}^m e_{k_l}(z_0,z_{i,1})}{q(z_0,z_{i,1})(z_{i,1}-z_0)}\Biggr)\Biggl( \prod_{i=1}^l\Biggl(\prod_{j=1}^{d_i-1}
 \frac{\prod_{l=1}^me_{k_l}(z_{i,j},z_{i,j+1})}{q(z_{i,j},z_{i,j+1})(2z_{i,j}\!-\!z_{i,j-1}\!-\!z_{i,j+1})\prod_{l=1}^m k_l z_{i,j}}
 \Biggr) \Biggr) .
\end{align*}
The integrand associated with the type (ii) loop graph $f_{\Gamma^{\mathrm{loop}}_d}$ is given by
\begin{align*}
&\frac{1}{2d(\prod_{i=1}^Na_i)^{d}} \left( \prod_{i=1}^d\frac{\prod_{l=1}^m e_{k_l}(z_{i},z_{i+1})}{q(z_i,z_{i+1})(2z_{i}-z_{i-1}-z_{i+1})\prod_{l=1}^m k_lz_{i}}\right) \left( \prod_{i=1}^d\frac{1}{(z_{i})^N} \right) \notag \\&\times\prod_{a=2}^{N-m-1}\left(\sum_{i=1}^d w_a(z_{i},z_{i+1})\right) ^{n_a}.
\end{align*}
The integrand associated with the type (iii) star graph of degree $d$ having a cluster vertex of degree $f\; (1\leq f \leq d-1)$ and the partition $\sigma=(d_1,d_2,\dots,d_l) \in P_{d-f}$ is given by
\begin{align*}
&\frac{\text{Sym}(\sigma)}{24\bigl(\prod_{i=1}^Na_i\bigr)^{d}(z_0)^{N(f-1)}} \Biggl( \prod_{i=1}^l\Biggl(\prod_{j=1}^{d_i}\frac{1}{(z_{i,j})^N}
 \Biggr) \Biggr)\biggl(-\frac{N-m}{N}\frac{1}{w^N}-\frac{N+m}{N}\frac{1}{(z_0)^N} \biggr)\notag \\
 &\qquad{}\times \frac{1}{(w-z_0)^2q(w,z_0)(q(z_0,z_0))^{f-1}} \Biggl( \prod_{p=1}^m\frac{1}{(k_pz_0)^{l-1}}\frac{e_{k_p}(w,z_0) (e_{k_p}(z_0,z_0) )^{f-1}}{k_p w (k_p z_0 )^{f}} \Biggl) \notag \\
 &\qquad{} \times \Biggl( \prod_{i=1}^l\frac{\prod_{p=1}^me_{k_p}(z_0,z_{i,1})}{q(z_0,z_{i,1})(z_{i,1}-z_0)}\Biggl)\\
 & \qquad{}\times \Biggl( \prod_{i=1}^l\Biggl(\prod_{j=1}^{d_i-1}\frac{\prod_{p=1}^m e_{k_p}(z_{i,j},z_{i,j+1})}{q(z_{i,j},z_{i,j+1})(2z_{i,j}-z_{i,j-1}-z_{i,j+1})\bigl(\prod_{p=1}^m k_pz_{i,j} \bigr)} \Biggr) \Biggr) \notag \\
&\times\prod_{a=2}^{N-m-1}\Biggl(w_a(w,z_0)+(f-1)w_a(z_0,z_0)+\sum_{i=1}^l \sum_{j=1}^{d_i}w_a(z_{i,j-1},z_{i,j})
 \Biggr) ^{n_a}.
\end{align*}
The integrand $f_{\Gamma^{\mathrm{point}}_d}$ is given by
\begin{align*}
&\frac{R(d)}{24\bigl(\prod_{i=1}^Na_i\bigr)^{d+1}(z_0)^{Nd+1}(q(z_0,z_0))^d} \Biggl( \prod_{l=1}^m \left( \frac{e_{k_l}(z_0,z_0)}{k_l z_0} \right)^d \Biggr)\Biggl( \prod_{a=2}^{N-m-1}(dw_a(z_0,z_0))^{n_a} \Biggr),
\end{align*}
where the symmetric factor $R(d)$ is defined as follows:
\begin{align*}
R(d)=\Biggl(\prod_{i=1}^Na_i\Biggr)\Biggl( \left(\frac{N-m}{2}\right) \frac{1}{d}-\Biggl( \sum_{j=1}^N \frac{1}{a_j}-\sum_{l=1}^m \frac{1}{k_l}\Biggr)\frac{1}{d^2}\Biggr).
\end{align*}
Lastly, we explain how to take the residue of $f_\Gamma$, according to our previous work~\cite{MK}.
\begin{Definition}
For each type of graph, the residue operation $\mathrm{Res}\colon f_\Gamma \to \mathbb{R}$ is defined as follows:
\begin{itemize}\itemsep=0pt
\item[(i)] First, we take the residue of $f_\Gamma$ at $z_0=0$. Next, we take the residue of the resulting function at $z_{i,j}=0$ and $z_{i,j}=\frac{z_{i,j+1}-z_{i,j-1}}{2}$, and sum them sequentially in ascending order of $j$ (for $1\leq j \leq d_{i}-1$). Lastly, we take the residue of the resulting function at $z_{i,d_i}=0$. The order among different $i$'s does not matter.

\item[(ii)] We take the residue of $f_\Gamma$ at $z_{i}=0$ and $z_{i}=\frac{z_{i+1}-z_{i-1}}{2}$, and sum them sequentially in ascending order of $i$, $1\leq i \leq d$.

\item[(iii)] First, we take the residue of $f_\Gamma$ at $w=z_0$. Then, the remaining process is the same as in the case of type~(i).

\item[(iv)] We take the residue of $f_\Gamma$ at $z_0=0$.
\end{itemize}
\end{Definition}

\section{Numerical tests for various examples} \label{Numerical Tests for Various Examples}
In this section, we present the numerical results for several $P(a_1,\dots,a_N\mid k_1,\dots,k_m)$'s obtained from Conjecture \ref{main}. Our examples consist of three types: degree $k$ hypersurfaces $P(a_1,a_2,\dots,a_N\mid k)$, complete intersections in the projective space, and general complete intersections $P(a_1,a_2,\dots,a_N\mid k_1,k_2,\dots,k_m)$. For simplicity, let us denote $e^{x^1}$ and $e^{t^1}$ by~$q$ and~$Q$, respectively.

\subsection{Examples of Fano hypersurfaces}

Below, we present the results for $P(1,1,1,2\mid 4)$, $P(1,1,1,1,2\mid 2)$, and $P(1,1,1,1,2\mid 4)$.
First, we show the results for $P(1,1,1,2\mid 4)$, a degree $4$ hypersurface in $P(1,1,1,2)$. This is a complex two-dimensional manifold. The mirror maps are given by
\begin{align*}
\begin{split}
&t^0=12 q+x^0+696 q^2 x^2+85344 q^3 \bigl(x^2\bigr)^2+\cdots, \\
&t^1=x^1+52 q x^2+4752 q^2 \bigl(x^2\bigr)^2+\frac{2193344}{3} q^3 \bigl(x^2\bigr)^3+\cdots, \\
&t^2=x^2+52 q \bigl(x^2\bigr)^2+\frac{22960}{3} q^2 \bigl(x^2\bigr)^3+1471808 q^3 \bigl(x^2\bigr)^4+\cdots.
\end{split}
\end{align*}
Their inversions, $x^0(t^*)$, $x^1(t^*)$, and $x^2(t^*)$, are given as follows:
\begin{align}
&x^0=t^0-12 Q-72 Q^2 t^2-864 Q^3 \bigl(t^2\bigr)^2+\cdots, \label{x^0k=4P(1,1,1,2)} \\
&x^1=t^1-52 Q t^2+656 Q^2 \bigl(t^2\bigr)^2-\frac{34720 Q^3 \bigl(t^2\bigr)^3}{3}+\cdots, \label{x^1k=4P(1,1,1,2)} \\
&x^2=t^2-52 Q \bigl(t^2\bigr)^2+\frac{1376 Q^2 \bigl(t^2\bigr)^3}{3}-\frac{167648 Q^3 \bigl(t^2\bigr)^4}{3}+\cdots. \label{x^2k=4P(1,1,1,2)}
\end{align}
The generating function $F_1^B$ is given by
\begin{align*}
F_1^B=-\frac{x^1}{12}-\frac{13}{3} q x^2-394 q^2 \bigl(x^2\bigr)^2-\frac{543920}{9} q^3 \bigl(x^2\bigr)^3+\cdots.
\end{align*}
Following Conjecture \ref{main}, we obtain the generating function $F_1^A$ by substituting (\ref{x^0k=4P(1,1,1,2)}), (\ref{x^1k=4P(1,1,1,2)}), and (\ref{x^2k=4P(1,1,1,2)}) into $F_{1}^{B}$:
\begin{align*}
F_1^A=-\frac{t^1}{12}+2 Q^2 \bigl(t^2\bigr)^2+\frac{224 Q^3 \bigl(t^2\bigr)^3}{3}+\cdots.
\end{align*}
Then, we obtain the genus $1$ Gromov--Witten invariants $\bigl\langle( {\mathcal O}_{h^2})^d\bigr\rangle_{1,d}$ as follows:
\begin{align*}
\langle {\mathcal O}_{h^2}\rangle_{1,1} =0, \qquad \bigl\langle ({\mathcal O}_{h^2})^2\bigr\rangle_{1,1}=2\times2!=4, \qquad \bigl\langle ({\mathcal O}_{h^2})^3\bigr\rangle_{1,1}=\frac{224}{3}\times3!=448.
\end{align*}
Next, we present the results for $P(1,1,1,1,2\mid 2)$, the degree $2$ hypersurface in $P(1,1,1,1,2)$. The mirror maps are given by
\begin{gather*}
t^0 =x^0+q\left(\frac{1}{6} \bigl(x^2\bigr)^3+ x^2x^3\right)\\
\hphantom{t^0 =}{} +q^2\left(\frac{32}{315} \bigl(x^2\bigr)^7+\frac{9}{5} \bigl(x^2\bigr)^5x^3+\frac{23}{3} \bigl(x^2\bigr)^3 \bigl(x^3\bigr)^2
 +\frac{13}{2} x^2\bigl(x^3\bigr)^3 \right)+\cdots,
\\
t^1 =x^1+q \left(\frac{\bigl(x^2\bigr)^4}{6}+\frac{3}{2} \bigl(x^2\bigr)^2x^3+\bigl(x^3\bigr)^2\right) \\
\hphantom{t^1 =}{}
 +q^2 \left(\frac{163 \bigl(x^2\bigr)^8}{1260}+\frac{121}{45} \bigl(x^2\bigr)^6x^3+\frac{89}{6} \bigl(x^2\bigr)^4\bigl(x^3\bigr)^2 +\frac{64}{3} \bigl(x^2\bigr)^2\bigl(x^3\bigr)^3 +\frac{11 \bigl(x^3\bigr)^4}{3}\right)+\cdots,
\\
t^2 =x^2+q \left(\frac{\bigl(x^2\bigr)^5}{12}+\frac{7}{6} \bigl(x^2\bigr)^3x^3+\frac{5}{2} x^2\bigl(x^3\bigr)^2 \right)\notag \\
\hphantom{t^2 =}{}
+q^2 \left(\frac{191 \bigl(x^2\bigr)^9}{2268}+\frac{193}{90} \bigl(x^2\bigr)^7x^3\!+16 \bigl(x^2\bigr)^5\bigl(x^3\bigr)^2\!+\frac{229}{6} \bigl(x^2\bigr)^3\bigl(x^3\bigr)^3+\frac{62}{3} x^2\bigl(x^3\bigr)^4\!\right) +\cdots,
\\
t^3 =x^3+q \left(\frac{\bigl(x^2\bigr)^6}{36}+\frac{7}{12} \bigl(x^2\bigr)^4x^3+\frac{5}{2} \bigl(x^2\bigr)^2\bigl(x^3\bigr)^2 +\frac{7}{6}\bigl(x^3\bigr)^3\right) \\
\hphantom{t^3 =}{}
+q^2 \left(\frac{106 \bigl(x^2\bigr)^{10}}{2835}+\frac{71}{60} \bigl(x^2\bigr)^8x^3+\frac{47}{4} \bigl(x^2\bigr)^6\bigl(x^3\bigr)^2 +\frac{1511}{36} \bigl(x^2\bigr)^4\bigl(x^3\bigr)^3+45 \bigl(x^2\bigr)^2\bigl(x^3\bigr)^4 \right. \\
\left.
\hphantom{t^3 =}{}
+\frac{19 \bigl(x^3\bigr)^5}{3}\right)+\cdots.
\end{gather*}
Inverting these mirror maps, we obtain
\begin{gather*}
x^0 =t^0+Q \left(-\frac{\bigl(t^2\bigr)^3}{6}-t^2t^3\right)\\
\hphantom{x^0 =}{} +Q^2 \left(-\frac{11 \bigl(t^2\bigr)^7}{2520}-\frac{2 \bigl(t^2\bigr)^5 t^3}{15}-\frac{13 \bigl(t^2\bigr)^3 \bigl(t^3\bigr)^2}{12}
-\frac{11 \bigl(t^2\bigr) \bigl(t^3\bigr)^3}{6}\right)+\cdots,
\\
x^1 =t^1+Q \left(-\frac{\bigl(t^2\bigr)^4}{6}-\frac{3 \bigl(t^2\bigr)^2 t^3}{2}-\bigl(t^3\bigr)^2\right)\notag \\
\hphantom{x^1 =}{}
 +Q^2 \left(-\frac{11 \bigl(t^2\bigr)^8}{2520}-\frac{83 \bigl(t^2\bigr)^6 t^3}{360}-\frac{13 \bigl(t^2\bigr)^4 \bigl(t^3\bigr)^2}{6}-\frac{49 \bigl(t^2\bigr)^2 \bigl(t^3\bigr)^3}{12}-\frac{\bigl(t^3\bigr)^4}{3}\right)+\cdots,
\\
x^2 =t^2+Q \left(-\frac{\bigl(t^2\bigr)^5}{12}-\frac{7 \bigl(t^2\bigr)^3 t^3}{6}-\frac{5 t^2 \bigl(t^3\bigr)^2}{2}\right) \\
\hphantom{x^2 =}{}
 +Q^2 \left(-\frac{29 \bigl(t^2\bigr)^9}{9072}-\frac{41 \bigl(t^2\bigr)^7 t^3}{180}-\frac{31 \bigl(t^2\bigr)^5 \bigl(t^3\bigr)^2}{12}-\frac{139 \bigl(t^2\bigr)^3 \bigl(t^3\bigr)^3}{18}-\frac{73 t^2 \bigl(t^3\bigr)^4}{12}\right)+\cdots,
\\
x^3 =t^3+Q \left(-\frac{\bigl(t^2\bigr)^6}{36}-\frac{7 \bigl(t^2\bigr)^4 t^3}{12}-\frac{5 \bigl(t^2\bigr)^2 \bigl(t^3\bigr)^2}{2}-\frac{7 \bigl(t^3\bigr)^3}{6}\right) \\
\hphantom{x^3 =}{}
 +Q^2 \left(-\frac{121 \bigl(t^2\bigr)^{10}}{45360}-\frac{127 \bigl(t^2\bigr)^8 t^3}{720}-\frac{173 \bigl(t^2\bigr)^6 \bigl(t^3\bigr)^2}{72}-\frac{95 \bigl(t^2\bigr)^4 \bigl(t^3\bigr)^3}{9}-\frac{41 \bigl(t^2\bigr)^2 \bigl(t^3\bigr)^4}{3}\right. \notag \\
 \left.\hphantom{x^3 =}{} -\frac{13 \bigl(t^3\bigr)^5}{12}\right)+\cdots.
\end{gather*}
The generating function $F_1^B$ is given by
\begin{gather*}
F_1^B =-\frac{x^1}{4}+q \left(-\frac{7}{144} \bigl(x^2\bigr)^4-\frac{5}{12} \bigl(x^2\bigr)^2x^3-\frac{7 \bigl(x^3\bigr)^2}{24}\right)+q^2 \left(-\frac{4541 \bigl(x^2\bigr)^8}{120960}-\frac{1097 \bigl(x^2\bigr)^6x^3 }{1440}\right. \\
\left.\hphantom{F_1^B =}{}
-\frac{33}{8} \bigl(x^2\bigr)^4\bigl(x^3\bigr)^2-\frac{853}{144} \bigl(x^2\bigr)^2\bigl(x^3\bigr)^3-\frac{19 \bigl(x^3\bigr)^4}{18}\right)+\cdots.
\end{gather*}
In the same way as the previous example, we obtain the generating function $F_1^A$ by substituting the inversions of the mirror maps into $F_{1}^{B}$,
\begin{gather*}
F_1^A =-\frac{t^1}{4}+Q \left(-\frac{\bigl(t^2\bigr)^4}{144}-\frac{\bigl(t^2\bigr)^2 t^3}{24}-\frac{\bigl(t^3\bigr)^2}{24}\right)\notag \\
\hphantom{F_1^A =}{}
+Q^2 \left(-\frac{23 \bigl(t^2\bigr)^8}{40320}-\frac{\bigl(t^2\bigr)^6 t^3}{160}-\frac{\bigl(t^2\bigr)^4 \bigl(t^3\bigr)^2}{48}-\frac{\bigl(t^2\bigr)^2 \bigl(t^3\bigr)^3}{48}\right)+\cdots.
\end{gather*}
We present numerical results of the genus $1$ Gromov--Witten invariants of this example in Table~\ref{P(1,1,1,1,2),k=2}. $N^0_{d,a,b}$, $N^1_{d,a,b}$, and $w_{a,b}$ denote $\bigl\langle ({\mathcal O}_{h^2})^a ({\mathcal O}_{h^3})^b\bigr\rangle_{0,d}$, $\bigl\langle ({\mathcal O}_{h^2})^a ({\mathcal O}_{h^3})^b\bigr\rangle_{1,d}$, and $w\bigl( ({\mathcal O}_{h^2})^a ({\mathcal O}_{h^3})^b\bigr)_{1,d}$, respectively. Since $P(1,1,1,1,k\mid k)$ is deformation
equivalent to $P(1,1,1,1)=CP^3$, these results agree with the results presented in~\cite{Getzler}.

\begin{table}[!ht]\small
\centering\renewcommand{\arraystretch}{1.29}
\caption{$P(1,1,1,1,2)$, $k=2$.} \label{P(1,1,1,1,2),k=2} \vspace{1mm}

\begin{tabular}{|c|r|r|r|r|r|}
\hline
$d$&\multicolumn{1}{|c|}{$(a,b)$} &\multicolumn{1}{|c|}{$N^0_{d,a,b}$} & \multicolumn{1}{|c|}{$N^1_{d,a,b}$} &\multicolumn{1}{|c|}{$\frac{2*d-1}{12}N^0_{d,a,b}+N^1_{d,a,b}$} &\multicolumn{1}{|c|}{$w_{a,b}$} \\
\hline
1&(0,2) & 1 &$-\frac{1}{12}$&0&$-\frac{7}{12}$\\ \hline
1&(2,1) & 1 &$-\frac{1}{12}$&0&$-\frac{5}{6}$\\ \hline
1&(4,0) & 2 &$-\frac{1}{6}$&0&$-\frac{7}{6}$\\ \hline
2&(0,4)& 0 & $0$& 0 &$-\frac{76}{3}$\\ \hline
2&(2,3)& 1 & $-\frac{1}{4}$& 0&$-\frac{853}{12}$ \\ \hline
2& (4,2)& 4 &$ -1$ &0&$-198$ \\ \hline
2& (6,1)& 18&$-\frac{9}{2}$& 0&$-\frac{1097}{2}$\\ \hline
2& (8,0) & 92&$-23$ & 0 &$-\frac{4541}{3}$ \\ \hline
3&(0,6)& 1 & $-\frac{5}{12}$& 0 &$-\frac{19959}{4}$\\ \hline
3&(2,5)& 5 & $-\frac{25}{12}$& 0&$-\frac{62338}{3}$ \\ \hline
3& (4,4)& 30&$ -\frac{25}{2}$&0&$-\frac{516827}{6}$ \\ \hline
3& (6,3)& 190&$-\frac{469}{6}$& 1&$-\frac{1068442}{3}$\\ \hline
3& (8,2) & 1312&$-\frac{1598}{3} $& 14&$-\frac{4408330}{3}$ \\ \hline
3& (10,1) & 9864&$-3960$ & 150 &$-\frac{18159922}{3}$\\ \hline
3& (12,0) & 80160&$-31900$ & 1500 &$-\frac{74719852}{3}$\\ \hline
4&(0,8)& 4 & $-\frac{4}{3} $&1&$-\frac{7111330}{3}$ \\ \hline
4& (2,7)&58&$-\frac{179}{6}$&4&$-\frac{26141813}{2}$ \\ \hline
4& (4,6)&480&$-248$& 32&$-71830274$\\ \hline
4&(6,5)& 4000 & $-\frac{6070}{3}$&310&$-\frac{1182256279}{3}$ \\ \hline
4& (8,4)& 35104&$-\frac{51772}{3}$&3220&$-2159333004$ \\ \hline
4& (10,3)&327888&$-156594 $& 34674&$-\frac{35458691818}{3}$\\ \hline
4& (12,2) & 3259680&$-1515824 $&385656&$-\frac{193936379144}{3}$ \\ \hline
4& (14,1) &34382544&$-15620216$ & 4436268 &$-353359995764$\\ \hline
4& (16,0) & 383306880&$-170763640 $ &52832040 &$-1930689790136$\\ \hline
5&(0,10)& 105 & $-\frac{147}{4} $&42&$-\frac{8363354113}{4}$ \\ \hline
5& (2,9)&1265&$-\frac{2379}{4}$&354&$-\frac{28682135389}{2}$ \\ \hline
5&(4,8)& 13354 & $-\frac{13047}{2} $&3492&$-\frac{196198477325}{2}$ \\ \hline
5& (6,7)&139098&$-\frac{132549}{2}$&38049&$-\frac{2010681907978}{3}$ \\ \hline
5& (8,6)&1492616&$-677808 $& 441654&$-\frac{13724961403006}{3}$ \\ \hline
5&(10,5)& 16744080 & $-7179606 $&5378454&$-\frac{93619004917238}{3}$ \\ \hline
5& (12,4)& 197240400&$-79637976 $&68292324&$-212735629674372$ \\ \hline
5& (14,3)&2440235712&$-928521900 $& 901654884&$-\frac{4348697671027760}{3}$\\ \hline
5& (16,2) & 31658432256&$-11385660384 $&12358163808&$-9873859605646752$ \\ \hline
5& (18,1) &429750191232&$-146713008096$ & 175599635328 &$-\frac{201722432909390752}{3}$\\ \hline
5& (20,0) & 6089786376960&$-1984020394752 $ &2583319387968&$-\frac{1373530281059327936}{3}$\\ \hline
\end{tabular}
\end{table}

Lastly, we consider the case of genus $1$ Gromov--Witten invariants for $P(1,1,1,1,2\mid 4)$, a~degree~$4$ hypersurface in $P(1,1,1,1,2)$. The mirror maps are given by
\begin{gather*}
t^0=x^0+12 q x^2+q^2 \bigl(848 \bigl(x^2\bigr)^3+2160 x^2x^3\bigr)+\cdots, \\
t^1=x^1+q \bigl(32 \bigl(x^2\bigr)^2+52 x^3\bigr)+q^2 \!\left(\frac{10568 \bigl(x^2\bigr)^4}{3}+16272 \bigl(x^2\bigr)^2x^3 +6416 \bigl(x^3\bigr)^2\right)+\cdots, \\
t^2=x^2+q \bigl(28 \bigl(x^2\bigr)^3+116 x^2x^3\bigl)+q^2\! \left(\!\frac{22936 \bigl(x^2\bigr)^5}{5}+\frac{96064}{3} \bigl(x^2\bigr)^3x^3+33552 x^2\bigl(x^3\bigr)^2 \!\right)\!+\cdots, \\
t^3=x^3+q \bigl(14 \bigl(x^2\bigr)^4+116 \bigl(x^2\bigr)^2x^3+84 \bigl(x^3\bigr)^2\bigl)\notag \\
\hphantom{t^3=}{} +q^2\! \left(\frac{17064 \bigl(x^2\bigr)^6}{5}+\frac{105904}{3} \bigl(x^2\bigr)^4x^3+74144 \bigl(x^2\bigr)^2\bigl(x^3\bigr)^2 +17808 \bigl(x^3\bigr)^3\right)+\cdots,
\end{gather*}
and their inversions are given as follows:
\begin{gather*}
x^0 =t^0-12 Q t^2+Q^2 \bigl(-128 \bigl(t^2\bigr)^3-144 t^2 t^3 \bigl)+\cdots, \\
x^1 =t^1+Q \bigl(-32 \bigl(t^2\bigr)^2-52 t^3\bigr)+Q^2 \left(\frac{64 \bigl(t^2\bigr)^4}{3}+512 \bigl(t^2\bigr)^2 t^3+656 \bigl(t^3\bigr)^2\right)+\cdots, \\
x^2 =t^2+Q \bigl(-28 \bigl(t^2\bigr)^3-116 t^2 t^3 \bigr)+Q^2 \left(\frac{1424 \bigl(t^2\bigr)^5}{5}-\frac{1216 \bigl(t^2\bigr)^3 t^3}{3}-4320 t^2 \bigl(t^3\bigr)^2\right)+\cdots, \\
x^3 =t^3+Q \bigl(-14 \bigl(t^2\bigr)^4-116 \bigl(t^2\bigr)^2 t^3-84 \bigl(t^3\bigr)^2\bigr) \\
\hphantom{x^3 =}{}
+Q^2 \left(\frac{1136 \bigl(t^2\bigr)^6}{5}-\frac{6184 \bigl(t^2\bigr)^4 t^3}{3}-9280 \bigl(t^2\bigr)^2 \bigl(t^3\bigr)^2+672 \bigl(t^3\bigr)^3\right)+\cdots.
\end{gather*}
Then, we present the generating functions $F_1^B$ and $F_1^A$ as follows:
\begin{gather*}
F_1^B =-\frac{x^1}{2}+q \bigl(-16 \bigl(x^2\bigr)^2-26 x^3\bigr) \\
\hphantom{F_1^B =}{}
+q^2 \left(-\frac{16096}{9} \bigl(x^2\bigr)^4-8184 \bigl(x^2\bigr)^2x^3-3212 \bigl(x^3\bigr)^2\right) +\cdots,
\end{gather*}
and
\begin{align*}
F_1^A =-\frac{t^1}{2}+Q^2 \left(-\frac{244 \bigl(t^2\bigr)^4}{9}-48 \bigl(t^2\bigr)^2 t^3-4 \bigl(t^3\bigr)^2\right)+\cdots.
\end{align*}
We present the genus $1$ Gromov--Witten invariants of this example in Table~\ref{P(1,1,1,1,2),k=4}.

\begin{table}[!ht]\small
\centering\renewcommand{\arraystretch}{1.29}
\centering
\caption{$P(1,1,1,1,2)$, $k=4$.}\label{P(1,1,1,1,2),k=4}\vspace{1mm}

\begin{tabular}{|c|r|r|r|r|r|}
\hline
$d$&\multicolumn{1}{|c|}{$(a,b)$} &\multicolumn{1}{|c|}{$N^0_{d,a,b}$} & \multicolumn{1}{|c|}{$N^1_{d,a,b}$} &\multicolumn{1}{|c|}{$\frac{d-1}{12}N^0_{d,a,b}+N^1_{d,a,b}$}&\multicolumn{1}{|c|}{$w_{a,b}$} \\
\hline
1&(0,1) & $24$ &0&0&$-26$\\ \hline
1&(2,0) & $80$ &0&0&$-32$\\ \hline
2&(0,2)&$144$ & $-8$& 4&$-6424$ \\ \hline
2& (2,1)&$1248$ &$-96$ &8&$-16368$ \\ \hline
2& (4,0)&$8192$&$-\frac{1952}{3} $& 32&$-\frac{128768}{3}$\\ \hline
3&(0,3)& $3456$ &$-128$& $448$&$-4177344$ \\ \hline
3& (2,2)&$48384$&$-5888$&$2176$&$-16223616$ \\ \hline
3& (4,1)&$491520$ &$-65792$& $16128$&$-\frac{193136768}{3}$\\ \hline
3& (6,0) & $5242880$&$-\frac{2206720}{3}$&$138240$&$-\frac{770396416}{3}$ \\ \hline
4&(0,4)& $165888 $ &$17664$& $59136$&$-4606798080$ \\ \hline
4&(2,3)& $3207168$ &$-279552 $& $522240$&$-24060080640$ \\ \hline
4& (4,2)&$44826624$&$-5371904 $&$5834752$&$-\frac{383428919296}{3}$ \\ \hline
4& (6,1)&$631504896$ &$-85794816 $& $72081408$&$-\frac{2046661990400}{3}$\\ \hline
4& (8,0) & $9330229248$&$-1381306368$&$951250944$&$-\frac{10951138650112}{3}$ \\ \hline
5&(0,5)& $12441600 $ &$5320704 $& $9467904$&$-7269250486272$ \\ \hline
5&(2,4)& $306892800 $ &$22683648$& $124981248$&$-47699671228416$ \\ \hline
5&(4,3)& $5506596864$ &$36335616$&$1871867904$&$-317367889719296$ \\ \hline
5& (6,2)&$97146372096$&$-2452013056 $&$29930110976$&$-\frac{6361733957066752}{3}$ \\ \hline
5& (8,1)&$1761381187584$ &$-82586042368 $& $504541020160$&$-\frac{42609768014790656}{3}$\\ \hline
5& (10,0) & $33262843985920$&$-\frac{6573683900416}{3}$&$8896386695168$&$-\frac{285858107179958272}{3}$ \\ \hline
\end{tabular}
\end{table}

\subsection{Example of Calabi--Yau threefolds}
In this subsection, we present the results for Calabi--Yau threefolds. When $\sum_{i=1}^{N}a_i=k$ (meaning the first Chern class is zero), the hypersurface $P(a_1,a_2,\dots,a_{N-1},a_N\mid k)$ becomes a Calabi--Yau manifold. The elliptic virtual structure constants for Calabi--Yau manifolds consist of contributions from type~(i),~(ii), and~(iv) graphs only, because the following proposition holds.

\begin{Proposition}[\cite{MK}]\label{Proposition1}
When the complete intersection is a Calabi--Yau manifold, for any posi\-tive $d$ and any $\Gamma \in \mathrm{Graph}^{\rm (iii)}_d$, $\mathrm{Res}(f_\Gamma)$ vanishes.
\end{Proposition}
The proof of this proposition will be given in Appendix~\ref{appendixA}.
In \cite{BCOV}, the genus~$1$ Gromov--Witten invariants for the Calabi--Yau 3-folds $P(1,1,1,1,2\mid 6)$, $P(1,1,1,1,4\mid 8)$, and $P(1,1,1,2,5\mid 10)$ were computed.
In the following, we present the numerical results for these three examples obtained using our formalism.
The mirror map and its inversion for the degree~$6$ hypersurface in $P(1,1,1,1,2)$ are given as follows:
\begin{gather*}
t=x+2772 q+9545850 q^2+53054643120 q^3+362147606012925 q^4+\cdots, \\
x=t-2772 Q-1861866 Q^2-5621359992 Q^3-20982861018549 Q^4+\cdots.
\end{gather*}
Since the mirror maps in Calabi--Yau manifolds involve only $t^1$ and $x^1$, we denote $t^1$ and $x^1$ by~$t$ and~$x$, respectively. The generating functions $F_{1}^{B}$ and $F_{1}^{A}$ are given by
\begin{gather*}
\begin{split}
& F_1^B=-\frac{7}{4} x-4194 q-14373450 q^2 - 80082321984 q^3 -547479376081866 q^4+\cdots, \\
& F_1^A=-\frac{7}{4}t+657 Q+\frac{1021167}{2} Q^2+1136816358 Q^3+\frac{18625762314603}{4}Q^4+\cdots.
\end{split}
\end{gather*}
The numerical results for the degree $8$ hypersurface in $P(1,1,1,1,4)$ are given as follows:
\begin{align*}
&t=x+15808 q+303422880 q^2+\frac{28300071331840 }{3}q^3+360758676442805200 q^4+\cdots, \\
&x=t-15808 Q - 53530016 Q^2-\frac{2907870121984}{3} Q^3 -20199602025147344 Q^4+\cdots,
\\
&F_1^B=-\frac{11 }{6}x-\frac{79568}{3} q-\frac{1519889680}{3} q^2-\frac{141843428301824}{9} q^3\notag \\
& \hphantom{F_1^B=}{} -\frac{1807425012871005968}{3} q^4+\cdots, \\
&F_1^A=-\frac{11}{6} t+\frac{7376}{3}Q+10778784 Q^2+\frac{1260969572864}{9} Q^3+3738046766828024 Q^4+\cdots.
\end{align*}
Lastly, we present the numerical results for the degree $10$ hypersurface in $P(1,1,1,2,5)$:
\begin{gather*}
t=x+179520 q+41513527200 q^2+15647390855936000 q^3 \\
\hphantom{t=}{}
+7272953267875497090000 q^4+\cdots, \\
x=t-179520 Q -9286096800 Q^2-1968068105216000 Q^3 \\
\hphantom{x=}{}-523504041681831810000 Q^4+\cdots, \\
F_1^B=-\frac{17 }{12}x-\frac{704320}{3} q-\frac{162228419200}{3} q^2-\frac{182678675797888000}{9} q^3 \\
\hphantom{F_1^B=}{}
-\frac{28070398497758257040000}{3} q^4+\cdots, \\
F_1^A=-\frac{17}{12} t+\frac{58640 }{3}Q+\frac{3677018600 }{3}Q^2+\frac{2727182857856000 }{9}Q^3 \\
\hphantom{F_1^A=}{}
+132217958645787677500 Q^4+\cdots.
\end{gather*}
At this stage, we briefly review how to compute the number of elliptic curves from these $F_{1}^{A}$.
In general, the perturbed two-point function $\langle{\mathcal O}_h {\mathcal O}_h\rangle(t)$ for a Calabi--Yau threefold given as a~hypersurface in weighted projective space has the structure
\begin{align*}
\langle{\mathcal O}_h {\mathcal O}_h\rangle(t)=\frac{kt}{\prod_{i=1}^N a_i}+\sum_{d=1}^\infty Q^d \langle{\mathcal O}_h {\mathcal O}_h\rangle_{0,d}.
\end{align*}
Let $n_{d}$ be the number of degree $d$ rational curves in the Calabi--Yau threefold. Then, $n_d$ is determined from the relation
 \begin{gather}
\frac{\mathrm{d}}{\mathrm{d}t}\langle{\mathcal O}_h {\mathcal O}_h\rangle(t)=\frac{k}{\prod_{i=1}^N a_i}+\sum_{d=1}^\infty \frac{n_d d^3 Q^d}{1-Q^d}. \label{number of rational curve}
\end{gather}
Let $m_d$ be the number of degree $d$ elliptic curves in the Calabi--Yau threefold. According to \cite{BCOV}, this $m_d$ is determined from the following relation:
 \begin{align}
F_1^A=-\frac{kc_{N-3}}{24\prod_{i=1}^N a_i}t-\sum_{d=1}^\infty m_d \log \Biggl( \prod_{p=1}^\infty \bigl(1-Q^{pd} \bigr)\Biggr) -\frac{1}{12} \sum_{d=1}^\infty n_d \log \bigl(1-Q^d \bigr). \label{BCOV formula}
\end{align}
We present the number of rational and elliptic curves for each case in Tables \ref{rationalcurveCalabi--Yauwaighted} and \ref{ellipticcurveCalabi--Yauwaighted}. Table~\ref{rationalcurveCalabi--Yauwaighted} is based on the results presented in \cite[p.~16]{AS}. The results for $m_{d}$, $d\leq 3$, agree with the results presented in~\cite[Table~2]{BCOV}.

\begin{table}[!ht]\small
\centering\renewcommand{\arraystretch}{1.29}
\centering
\caption{Rational curves on Calabi--Yau weighted projective space.}\label{rationalcurveCalabi--Yauwaighted}\vspace{1mm}

\begin{tabular}{|c|r|r|r|r|r|}
\hline
$d$&\multicolumn{1}{|c|}{$P(1,1,1,1,2)$, $k=6$} & \multicolumn{1}{|c|}{$P(1,1,1,1,4)$, $k=8$} &\multicolumn{1}{|c|}{$P(1,1,1,2,5)$, $k=10$} \\
\hline
$1$&$7884$ & $29504$ &$231200$\\ \hline
$2$&$6028452$& $128834912$ &$12215785600$\\ \hline
$3$&$11900417220$&$1423720546880$ & $1700894366474400$ \\ \hline
$4$&$34600752005688$ & $23193056024793312$ &$350154658851324656000$\\ \hline
$5$&$124595034333130080$& $467876474625249316800$ &$89338191421813572850115680$\\ \hline
\end{tabular}
\end{table}

\begin{table}[!ht]\small
\centering\renewcommand{\arraystretch}{1.2}
\centering
\caption{Elliptic curves on Calabi--Yau weighted projective space.}\label{ellipticcurveCalabi--Yauwaighted}\vspace{1mm}

\begin{tabular}{|c|r|r|r|}
\hline
$d$&\multicolumn{1}{|c|}{$P(1,1,1,1,2)$, $k=6$} & \multicolumn{1}{|c|}{$P(1,1,1,1,4)$, $k=8$} &\multicolumn{1}{|c|}{$P(1,1,1,2,5)$, $k=10$} \\
\hline
$1$&$0$ & $0$ &$280$\\ \hline
$2$&$7884$& $41312$ &$207680680$\\ \hline
$3$&$145114704$&$21464350592$ & $161279120326560$ \\ \hline
$4$&$1773044315001$ & $1805292092664544$ &$103038403740690105440$\\ \hline
$5$&$17144900584158168$& $101424054914016355712$ &$59221844124053623534386928$\\ \hline
\end{tabular}
\end{table}

\subsection{Examples of complete intersections in projective space}
In this subsection, we present the results for complete intersections in projective space ${CP}^{N-1}$, denoted as $P(1,1,\dots,1\mid k_1,k_2,\dots,k_m)$. Since $a_1=a_2=\cdots=a_N=1$, we abbreviate $P(1,1,\dots,1\mid k_1,k_2,\dots,k_m)$ as $(k_1,k_2,\dots,k_m)_N$. We now introduce four examples: $(2,2)_5$, $(2,2,2)_6$, $(2,2,2)_7$, and $(2,2,3)_7$.

First, we present the results for $(2,2)_5$.
The mirror maps are given by
\begin{gather*}
t^0=x^0+4 q+56 q^2 x^2+1696 q^3 \bigl(x^2\bigr)^2+\cdots, \\
t^1=x^1+12 q x^2+272 q^2 \bigl(x^2\bigr)^2+10432 q^3 \bigl(x^2\bigr)^3+\cdots, \\
t^2=x^2+12 q \bigl(x^2\bigr)^2+\frac{1328}{3} q^2 \bigl(x^2\bigr)^3+\frac{64064}{3} q^3 \bigl(x^2\bigr)^4+\cdots.
\end{gather*}
Their inversions are given as follows:
\begin{gather*}
x^0=t^0-4 Q-8 Q^2 t^2-32 Q^3 \bigl(t^2\bigr)^2+\cdots, \\
x^1=t^1-12 Q t^2+16 Q^2 \bigl(t^2\bigr)^2-32 Q^3 \bigl(t^2\bigr)^3+\cdots, \\
x^2=t^2-12 Q \bigl(t^2\bigr)^2-\frac{32 Q^2 \bigl(t^2\bigr)^3}{3}-\frac{2336 Q^3 \bigl(t^2\bigr)^4}{3}+\cdots.
\end{gather*}
The generating functions are given by
\begin{gather*}
F_1^B=-\frac{x^1}{6}-2 q x^2-\frac{136}{3} q^2 \bigl(x^2\bigr)^2-\frac{5216}{3} q^3 \bigl(x^2\bigr)^3+\cdots, \qquad
F_1^A=-\frac{t^1}{6}+\cdots.
\end{gather*}
We present the genus $1$ Gromov--Witten invariants obtained from these results in Table \ref{(2,2)CP4}, where~$w$ is the value of the corresponding elliptic virtual structure constant.

\begin{table}[!ht]\small
\centering\renewcommand{\arraystretch}{1.2}
\centering
\caption{$(2,2)_5$.}\label{(2,2)CP4}\vspace{1mm}

\begin{tabular}{|c|r|r|r|}
\hline
$d$&\multicolumn{1}{|c|}{$N^0_{d}$} & \multicolumn{1}{|c|}{$N^1_{d}$} &\multicolumn{1}{|c|}{$w$} \\
\hline
1 & $16$ &0&$-2$\\ \hline
2&$40$ & $0$&$-\frac{272}{3}$ \\ \hline
3& $256$ &$0$&$-10432$ \\ \hline
4& $3328$ &$256$&$-\frac{6007040}{3}$ \\ \hline
5& $69632$ &$16384$& $-\frac{1633808384}{3}$ \\ \hline
\end{tabular}
\end{table}

Next, we present the results for $(2,2,2)_6$, which is a K3 surface.
The mirror map is given by
\begin{align*}
t&=x+24 q + 564 q^2 + 19904 q^3 +\cdots, \\
x&=t-24 Q + 12 Q^2 - 32 Q^3 +\cdots.
\end{align*}
All elliptic virtual structure constants up to degree $5$ vanish. Therefore, the generating functions up to degree $5$ are given by $F_1^B=0$, $F_1^A=0$.
This result agrees with the discussions in \cite{BCOV}, where the vanishing of genus $1$ Gromov--Witten invariants for the K3 surface is suggested.

We then turn to the results for $(2,2,2)_7$, which is a Fano threefold. We present the results up to degree $3$.
The mirror maps are given by
\begin{gather*}
t^0=x^0+8 q+368 q^2 x^2+q^3 \bigl(40320 \bigl(x^2\bigr)^2+22656 x^3\bigl)+\cdots, \\
t^1=x^1+32 q x^2+q^2 \bigl(2560 \bigl(x^2\bigr)^2+1936 x^3\bigl)+q^3 \left(\frac{1113344 \bigl(x^2\bigr)^3}{3}+561920 x^2x^3\right)+\cdots, \\
t^2=x^2+q \bigl(40 \bigl(x^2\bigr)^2+56 x^3\bigl)+q^2 \left(\frac{16960 \bigl(x^2\bigr)^3}{3}+11936 x^2 x^3\right) \\
\hphantom{t^2=}{} +q^3 \bigl(1073152 \bigl(x^2\bigr)^4+3028224 \bigl(x^2\bigr)^2
x^3+714496 \bigl(x^3\bigr)^2\bigr)+\cdots, \\
t^3=x^3+q \left(\frac{80 \bigl(x^2\bigr)^3}{3}+112 x^2 x^3\right)+q^2 \bigl(6400 \bigl(x^2\bigr)^4+26560 \bigl(x^2\bigr)^2 x^3+9248 \bigl(x^3\bigr)^2\bigl) \\
\hphantom{t^3=}{}
+q^3 \left(\frac{7816192
\bigl(x^2\bigr)^5}{5}+\frac{21595136}{3} \bigl(x^2\bigr)^3 x^3+5016576 x^2 \bigl(x^3\bigr)^2\right)+\cdots.
\end{gather*}
Their inversions are given as follows:
\begin{gather*}
x^0=t^0-8 Q-112 Q^2 t^2+Q^3 \bigl(-4096 \bigl(t^2\bigr)^2-896 t^3\bigr)+\cdots, \\
x^1=t^1-32 Q t^2+Q^2 \bigl(-256 \bigl(t^2\bigr)^2-144 t^3\bigl)+Q^3 \bigl(-10240 \bigl(t^2\bigr)^3-6656 t^2 t^3\bigr)+\cdots, \\
x^2=t^2+Q \bigl(-40 \bigl(t^2\bigr)^2-56 t^3\bigr)+Q^2 \bigl(320 \bigl(t^2\bigr)^3+608 t^2 t^3\bigr) \\
\hphantom{x^2=}{}+Q^3 \bigl(-512 \bigl(t^2\bigr)^4+1664
\bigl(t^2\bigr)^2 t^3+3200 \bigl(t^3\bigr)^2\bigr)+\cdots, \\
x^3=t^3+Q \left(-\frac{80 \bigl(t^2\bigr)^3}{3}-112 t^2 t^3\right)+Q^2 \left(640 \bigl(t^2\bigr)^4-1472 \bigl(t^2\bigr)^2 t^3-2976 \bigl(t^3\bigr)^2\right) \\
\hphantom{x^3=}{}+Q^3
\left(\frac{145408 \bigl(t^2\bigr)^5}{5}-\frac{396544 \bigl(t^2\bigr)^3 t^3}{3}-251136 t^2 \bigl(t^3\bigr)^2\right)+\cdots.
\end{gather*}
The generating functions are given by
\begin{gather*}
F_1^B =-x^1-\frac{80}{3} q x^2+q^2 \left(-2176 \bigl(x^2\bigr)^2-\frac{4912 x^3}{3}\right)\notag \\
\hphantom{F_1^B =}{}
+q^3 \left(-\frac{958720}{3} \bigl(x^2\bigr)^3-479488 x^2 x^3\right)+\cdots, \\
F_1^A =-t^1+\frac{16 t^2 Q}{3}+\left(-\frac{16384 \bigl(t^2\bigr)^3}{9}-\frac{3328 t^2 t^3}{3}\right) Q^3+\cdots.
\end{gather*}
We present the results for genus $1$ Gromov--Witten invariants in Table~\ref{(2,2,2)CP6}, where $N^0_{d,a,b}$ represents $\bigl\langle({\mathcal{O}}_{h^2})^a({\mathcal{O}}_{h^3})^b\bigr\rangle_{0,d}$
and $N^1_{d,a,b}$ represents $\bigl\langle({\mathcal{O}}_{h^2})^a({\mathcal{O}}_{h^3})^b\bigr\rangle_{1,d}$.

\begin{table}[!ht]\small
\centering\renewcommand{\arraystretch}{1.29}
\centering
\caption{$(2,2,2)_7$.}\label{(2,2,2)CP6}\vspace{1mm}

\begin{tabular}{|c|r|r|r|r|r|}
\hline
$d$&\multicolumn{1}{|c|}{$(a,b)$}& \multicolumn{1}{|c|}{$N^0_{d,a,b}$} & \multicolumn{1}{|c|}{$N^1_{d,a,b}$} &\multicolumn{1}{|c|}{$\frac{d-2}{24}N^0_{d,a,b}+N^1_{d,a,b}$}& \multicolumn{1}{|c|}{$w_{a,b}$} \\
\hline
1&(1,0) & $128$ &$\frac{16}{3}$&0&$-\frac{80}{3}$\\ \hline
2&(0,1)&$608$ & $0$&0&$-\frac{4912}{3}$ \\ \hline
2& (2,0)&$3200$ &$0$ &0&$-4352$ \\ \hline
3&(1,1)& $26624$ &$-\frac{3328 }{3}$& $0$&$-479488$ \\ \hline
3& (3,0)&$262144$&$-\frac{32768}{3}$&$0$&$-1917440$ \\ \hline
4&(0,2)& $242176$ &$-\frac{57856}{3}$& $896$&$-\frac{150728192}{3}$ \\ \hline
4&(2,1)& $2914304$ &$-\frac{696320 }{3} $& $10752$&$-\frac{808314880}{3}$ \\ \hline
4& (4,0)&$41943040$&$-\frac{10141696}{3}$&$114688$&$-\frac{4322443264}{3}$ \\ \hline
5&(1,2)& $33062912$ &$-3444736 $& $688128$&$-\frac{103377682432}{3}$ \\ \hline
5&(3,1)& $549453824$ &$-57344000 $& $11337728$&$-\frac{693950070784}{3}$ \\ \hline
5&(5,0)& $10401873920$ &$-1118306304$&$181927936$&$-\frac{4649625714688}{3}$ \\ \hline
\end{tabular}
\end{table}

Lastly, we present the result for $(2,2,3)_7$, which is a Calabi--Yau threefold.
The mirror map and its inversion are given by
\begin{gather*}
t=x+108 q + 14202 q^2 + 2974032 q^3 +\cdots, \\
x=t-108 Q - 2538 Q^2 - 262152 Q^3 +\cdots.
\end{gather*}
The generating functions are given as follows:
\begin{align*}
\begin{split}
& F_1^B=-\frac{5 x}{2}-210 q - 27126 q^2 - 5688480 q^3+\cdots, \\
& F_1^A=-\frac{5 t}{2} + 60 Q + 1899 Q^2 + 134376 Q^3 +\cdots.
\end{split}
\end{align*}
From this $F_{1}^{A}$, we obtain the number of rational curves $n_d$ and the number of elliptic curves $m_d$ by using equations~(\ref{number of rational curve}) and~(\ref{BCOV formula}).
The results are presented in Table~\ref{(2,2,3)CP6}.

\begin{table}[!ht]\small\vspace{-2mm}
\centering\renewcommand{\arraystretch}{1.2}
\centering
\caption{$n_d$ and $m_d$ in degrees $(2,2,3)_7$.}\label{(2,2,3)CP6}\vspace{1mm}

\begin{tabular}{|c|r|r|r|r|r|}
\hline
$d$&\multicolumn{1}{|c|}{$n_d$} & \multicolumn{1}{|c|}{$m_d$} \\
\hline
1& $720$ &$0$\\ \hline
2&$22428$ & $0$ \\ \hline
3& $1611504$ &$64$ \\ \hline
4& $168199200$ &$265113$ \\ \hline
5& $21676931712$ &$198087264$ \\ \hline
\end{tabular}\vspace{-2mm}
\end{table}

\subsection{Examples of complete intersections in weighted projective space}
 In this subsection, we present the results for complete intersections in weighted projective space. Specifically, we compute the two examples:
 $P(1,1,1,1,2\mid 2,2)$ and $P(1,1,1,1,1,2\mid 2,2)$.
The results for $P(1,1,1,1,2\mid 2,2)$ are summarized in Table~\ref{(1,1,1,1,2|2,2)}, where $n_d=\bigl\langle ({\mathcal {O}}_{h^2})^{2d-1}\bigr\rangle_{0,d}$ and $m_d=\bigl\langle ({\mathcal {O}}_{h^2})^{2d}\bigr\rangle_{1,d}$. Since the dimension of $P(1,1,1,1,2\mid 2,2)$ is two, the mirror maps are given~by
\begin{gather*}
t^0=x^0+2 q x^2+10 q^2 \bigl(x^2\bigr)^3+\frac{320}{3} q^3 \bigl(x^2\bigr)^5+\cdots,\\
t^1=x^1+3 q \bigl(x^2\bigr)^2+\frac{131}{6} q^2 \bigl(x^2\bigr)^4+\frac{12329}{45} q^3 \bigl(x^2\bigr)^6+\cdots, \\
t^2=x^2+2 q \bigl(x^2\bigr)^3+\frac{313}{15} q^2 \bigl(x^2\bigr)^5+\frac{10764}{35} q^3 \bigl(x^2\bigr)^7+\cdots,
\end{gather*}
and their inversions are
\begin{gather*}
\begin{split}
&x^0=t^0-2 Qt^2-\frac{4 Q^3 \bigl(t^2\bigr)^5}{15}+\cdots,\\
&x^1=t^1-3 Q \bigl(t^2\bigr)^2-\frac{5 Q^2 \bigl(t^2\bigr)^4}{6}-\frac{91 Q^3 \bigl(t^2\bigr)^6}{9}+\cdots, \\
& x^2=t^2-2 Q \bigl(t^2\bigr)^3-\frac{43 Q^2 \bigl(t^2\bigr)^5}{15}-\frac{500 Q^3 \bigl(t^2\bigr)^7}{21}+\cdots.
\end{split}
\end{gather*}
The generating functions are given by
\begin{gather*}
F_1^B=-\frac{x^1}{6}-\frac{1}{2} q \bigl(x^2\bigr)^2-\frac{131}{36} q^2 \bigl(x^2\bigr)^4-\frac{12329}{270} q^3 \bigl(x^2\bigr)^6+\cdots,\qquad
F_1^A=-\frac{t^1}{6}+\cdots.
\end{gather*}
Since $P(1,1,1,1,2\mid 2,2)$ is biholomorphically equivalent to a degree~$2$ hypersurface in ${CP}^{3}$, the result agrees with the corresponding result in~\cite{MK}.

\begin{table}[!ht]\small
\centering\renewcommand{\arraystretch}{1.2}
\centering
\caption{$P(1,1,1,1,2\mid 2,2)$.}\label{(1,1,1,1,2|2,2)}\vspace{1mm}
\begin{tabular}{|c|r|r|r|}
\hline
$d$&\multicolumn{1}{|c|}{$n_d$} & \multicolumn{1}{|c|}{$m_d$} \\
\hline
1& $4$ &$0$\\ \hline
2&$8$ & $0$ \\ \hline
3& $64 $ &$0$ \\ \hline
4& $1792$ &$256 $ \\ \hline
5& $99328 $ &$40960$ \\ \hline
\end{tabular}
\end{table}

Then, we present the results for $P(1,1,1,1,1,2\mid 2,2)$, which is a Fano threefold. The results for the Gromov--Witten invariants of $P(1,1,1,1,1,2\mid 2,2)$ are presented in Table~\ref{(1,1,1,1,1,2|2,2)}. In this table, $N^0_{d,a,b}=\bigl\langle (\mathcal{O}_{h^2})^{a} (\mathcal{O}_{h^3})^{b}\bigr\rangle_{0,d}$ and $m_d=\bigl\langle (\mathcal{O}_{h^2})^{a} (\mathcal{O}_{h^3})^{b}\bigr\rangle_{1,d}$, respectively.
The selection rules are given by
\begin{gather*}
\bigl\langle (\mathcal{O}_{h^2})^{a} (\mathcal{O}_{h^3})^{b}\bigr\rangle_{0,d} \neq0 \ \Longrightarrow \ 3d=a+2b, \\
\bigl\langle (\mathcal{O}_{h^2})^{a} (\mathcal{O}_{h^3})^{b}\bigr\rangle_{1,d} \neq0 \ \Longrightarrow \ 3d=a+2b.
\end{gather*}
The mirror maps are given by
\begin{gather*}
t^0 =x^0+q \bigl(\bigl(x^2\bigr)^2+2 x^3\bigr)+q^2 \left(\frac{58 \bigl(x^2\bigr)^5}{15}+\frac{98}{3} \bigl(x^2\bigr)^3 x^3+42 x^2 \bigl(x^3\bigr)^2\right) \\
\hphantom{t^0 =}{}
 +q^3 \left(\frac{4201 \bigl(x^2\bigr)^8}{126}+\frac{7546}{15}
\bigl(x^2\bigr)^6 \bigl(x^3\bigr)+\frac{6110}{3} \bigl(x^2\bigr)^4 \bigl(x^3\bigr)^2\right.\\
\left.\hphantom{t^0 =}{} +\frac{6634}{3} \bigl(x^2\bigr)^2 \bigl(x^3\bigr)^3+301 \bigl(x^3\bigr)^4\right)+\cdots,
\\
t^1 =x^1+q \left(\frac{4 \bigl(x^2\bigr)^3}{3}+6 x^2 x^3\right)\\
\hphantom{t^1 =}{} +q^2 \left(\frac{326 \bigl(x^2\bigr)^6}{45}+\frac{244}{3} \bigl(x^2\bigr)^4 x^3+182 \bigl(x^2\bigr)^2 \bigl(x^3\bigr)^2+\frac{134
\bigl(x^3\bigr)^3}{3}\right) \\
\hphantom{t^1 =}{} +q^3 \left(\frac{13499 \bigl(x^2\bigr)^9}{189}+\frac{26812}{21} \bigl(x^2\bigr)^7 x^3+\frac{33252}{5} \bigl(x^2\bigr)^5 \bigl(x^3\bigr)^2+\frac{98872}{9} \bigl(x^2\bigr)^3 \bigl(x^3\bigr)^3\right. \notag \\
\left.\hphantom{t^1 =}{} +\frac{12196}{3} x^2 \bigl(x^3\bigr)^4\right)+\cdots,
\\
t^2 =x^2+q \left(\frac{5 \bigl(x^2\bigr)^4}{6}+7 \bigl(x^2\bigr)^2 x^3+5 \bigl(x^3\bigr)^2\right)\notag \\
\hphantom{t^2 =}{} +q^2 \left(\frac{94 \bigl(x^2\bigr)^7}{15}+\frac{1393}{15} \bigl(x^2\bigr)^5 x^3+\frac{986}{3}
\bigl(x^2\bigr)^3 \bigl(x^3\bigr)^2+234 x^2 \bigl(x^3\bigr)^3\right) \notag \\
\hphantom{t^2 =}{} +q^3 \left(\frac{666611 \bigl(x^2\bigr)^{10}}{9450}+\frac{21069}{14} \bigl(x^2\bigr)^8 x^3+\frac{450524}{45} \bigl(x^2\bigr)^6 \bigl(x^3\bigr)^2+\frac{214577}{9}
\bigl(x^2\bigr)^4 \bigl(x^3\bigr)^3\right. \notag \\
 \left.
\hphantom{t^2 =}{}
 +17085 \bigl(x^2\bigr)^2 \bigl(x^3\bigr)^4+\frac{24577 \bigl(x^3\bigr)^5}{15}\right)+\cdots,
\\
t^3 =x^3+q \left(\frac{\bigl(x^2\bigr)^5}{3}+\frac{14}{3} \bigl(x^2\bigr)^3 x^3+10 x^2 \bigl(x^3\bigr)^2\right)\notag \\
\hphantom{t^3 =}{} +q^2 \left(\frac{629 \bigl(x^2\bigr)^8}{180}+\frac{1027}{15} \bigl(x^2\bigr)^6 x^3+\frac{1088}{3}
\bigl(x^2\bigr)^4 \bigl(x^3\bigr)^2+518 \bigl(x^2\bigr)^2 \bigl(x^3\bigr)^3+\frac{280 \bigl(x^3\bigr)^4}{3}\right)\notag \\
\hphantom{t^3 =}{} +q^3 \left(\frac{2362258 \bigl(x^2\bigr)^{11}}{51975}+\frac{41093}{35} \bigl(x^2\bigr)^9 x^3+\frac{3145468}{315}
\bigl(x^2\bigr)^7 \bigl(x^3\bigr)^2+\frac{1495106}{45} \bigl(x^2\bigr)^5 \bigl(x^3\bigr)^3\right.\notag \\
 \left.
 \hphantom{t^3 =}{}
 +39738 \bigl(x^2\bigr)^3 \bigl(x^3\bigr)^4+\frac{57354}{5} x^2 \bigl(x^3\bigr)^5\right)+\cdots,
\end{gather*}
and their inversions are given as follows:
\begin{gather*}
x^0 =t^0+Q \left(-\bigl(t^2\bigr)^2-2 t^3\right)+Q^2 \left(-\frac{\bigl(t^2\bigr)^5}{5}-\frac{2 \bigl(t^2\bigr)^3 t^3}{3}\right)\notag \\
\hphantom{x^0 =}{}
 +Q^3 \left(-\frac{367
\bigl(t^2\bigr)^8}{1260}-\frac{31 \bigl(t^2\bigr)^6 t^3}{15}-\frac{13 \bigl(t^2\bigr)^4 \bigl(t^3\bigr)^2}{3}-\frac{8 \bigl(t^2\bigr)^2 \bigl(t^3\bigr)^3}{3}\right)+\cdots,
\\
x^1 =t^1+Q \left(-\frac{4 \bigl(t^2\bigr)^3}{3}-6 t^2 t^3\right)+Q^2 \left(-\frac{2 \bigl(t^2\bigr)^6}{15}-\frac{13 \bigl(t^2\bigr)^4 t^3}{3}-24
\bigl(t^2\bigr)^2 \bigl(t^3\bigr)^2-\frac{44 \bigl(t^3\bigr)^3}{3}\right)\notag \\
\hphantom{x^1 =}{} +Q^3 \left(-\frac{989 \bigl(t^2\bigr)^9}{1890}-\frac{778 \bigl(t^2\bigr)^7 t^3}{35}-\frac{3338 \bigl(t^2\bigr)^5
\bigl(t^3\bigr)^2}{15}-\frac{5968 \bigl(t^2\bigr)^3 \bigl(t^3\bigr)^3}{9}\right.\\
\left.\hphantom{x^1 =}{} -\frac{1192 t^2 \bigl(t^3\bigr)^4}{3}\right)+\cdots,
\\
x^2 =t^2+Q \left(-\frac{5 \bigl(t^2\bigr)^4}{6}-7 \bigl(t^2\bigr)^2 t^3-5 \bigl(t^3\bigr)^2\right)\notag \\
\hphantom{x^2 =}{}
+Q^2 \left(-\frac{2 \bigl(t^2\bigr)^7}{45}-\frac{113
\bigl(t^2\bigr)^5 t^3}{15}-\frac{146 \bigl(t^2\bigr)^3 \bigl(t^3\bigr)^2}{3}-34 t^2 \bigl(t^3\bigr)^3\right)\notag \\
\hphantom{x^2 =}{} +Q^3 \left(-\frac{2497 \bigl(t^2\bigr)^{10}}{3150}-\frac{12365
\bigl(t^2\bigr)^8 t^3}{252}-547 \bigl(t^2\bigr)^6 \bigl(t^3\bigr)^2-1822 \bigl(t^2\bigr)^4 \bigl(t^3\bigr)^3\right. \notag \\
\left.
\hphantom{x^2 =}{}-1432 \bigl(t^2\bigr)^2 \bigl(t^3\bigr)^4
-\frac{684 \bigl(t^3\bigr)^5}{5}\right)+\cdots,
\\
x^3 =t^3+Q \left(-\frac{\bigl(t^2\bigr)^5}{3}-\frac{14 \bigl(t^2\bigr)^3 t^3}{3}-10 t^2 \bigl(t^3\bigr)^2\right)\notag \\
\hphantom{x^3 =}{} +Q^2 \left(-\frac{19 \bigl(t^2\bigr)^8}{180}-\frac{127
\bigl(t^2\bigr)^6 t^3}{15}-\frac{200 \bigl(t^2\bigr)^4 \bigl(t^3\bigr)^2}{3}-118 \bigl(t^2\bigr)^2 \bigl(t^3\bigr)^3-\frac{130 \bigl(t^3\bigr)^4}{3}\right) \notag \\
\hphantom{x^3 =}{} +Q^3 \left(-\frac{48793
\bigl(t^2\bigr)^{11}}{51975}-\frac{22417 \bigl(t^2\bigr)^9 t^3}{378}-\frac{236618 \bigl(t^2\bigr)^7 \bigl(t^3\bigr)^2}{315}-\frac{144296 \bigl(t^2\bigr)^5 \bigl(t^3\bigr)^3}{45} \right. \notag \\
 \left.\hphantom{x^3 =}{} -\frac{41120
\bigl(t^2\bigr)^3 \bigl(t^3\bigr)^4}{9}-\frac{26012 t^2 \bigl(t^3\bigr)^5}{15}\right)+\cdots.
\end{gather*}
The generating functions are the following:
\begin{gather*}
F_1^B =-\frac{x^1}{3}+q \left(-\frac{1}{2} \bigl(x^2\bigr)^3-\frac{13}{6} x^2 x^3\right)\notag \\
\hphantom{F_1^B =}{}
+q^2 \left(-\frac{1489}{540} \bigl(x^2\bigr)^6-\frac{1087}{36} \bigl(x^2\bigr)^4 x^3-66
\bigl(x^2\bigr)^2 \bigl(x^3\bigr)^2-\frac{287 \bigl(x^3\bigr)^3}{18}\right)  \\
\hphantom{F_1^B =}{}+q^3 \left(-\frac{461369 \bigl(x^2\bigr)^9}{17010}-\frac{450599}{945} \bigl(x^2\bigr)^7 x^3-\frac{109972}{45} \bigl(x^2\bigr)^5 \bigl(x^3\bigr)^2-\frac{107299}{27}
\bigl(x^2\bigr)^3 \bigl(x^3\bigr)^3\right. \notag \\
 \left.
\hphantom{F_1^B =}{}
 -\frac{26125}{18} x^2 \bigl(x^3\bigr)^4\right)+\cdots,
\\
F_1^A =-\frac{t^1}{3}+Q \left(-\frac{\bigl(t^2\bigr)^3}{18}-\frac{t^2 t^3}{6}\right)+Q^2 \left(-\frac{2 \bigl(t^2\bigr)^6}{27}-\frac{4
\bigl(t^2\bigr)^4 t^3}{9}-\frac{2 \bigl(t^2\bigr)^2 \bigl(t^3\bigr)^2}{3}-\frac{2 \bigl(t^3\bigr)^3}{9}\right)\notag \\
\hphantom{F_1^A =}{}
 +Q^3 \left(-\frac{121 \bigl(t^2\bigr)^9}{1215}-\frac{118 \bigl(t^2\bigr)^7
t^3}{135}-\frac{112 \bigl(t^2\bigr)^5 \bigl(t^3\bigr)^2}{45}-\frac{70 \bigl(t^2\bigr)^3 \bigl(t^3\bigr)^3}{27}-\frac{7 t^2 \bigl(t^3\bigr)^4}{9}\right)+\cdots.
\end{gather*}
Since $P(1,1,1,1,1,2\mid 2,2)$ is deformation equivalent to a degree $2$ hypersurface in ${CP}^{4}$, the result agrees with the corresponding result in~\cite{MK}.

\begin{table}[!ht]\small
\centering\renewcommand{\arraystretch}{1.29}
\centering
\caption{$P(1,1,1,1,1,2\mid 2,2)$.}\label{(1,1,1,1,1,2|2,2)}\vspace{1mm}

\begin{tabular}{|c|r|r|r|r|r|}
\hline
$d$&\multicolumn{1}{|c|}{$(a,b)$}& \multicolumn{1}{|c|}{$N^0_{d,a,b}$} & \multicolumn{1}{|c|}{$N^1_{d,a,b}$} &\multicolumn{1}{|c|}{$\frac{3d-2}{24}N^0_{d,a,b}+N^1_{d,a,b}$}& \multicolumn{1}{|c|}{$w_{a,b}$} \\
\hline
1&(1,1) & $4$ &$-\frac{1}{6}$&0&$-\frac{13}{6}$\\ \hline
1&(3,0) & $8$ &$-\frac{1}{3}$&0&$-3$\\ \hline
2&(0,3)&$8$ & $-\frac{4}{3}$&0&$-\frac{287}{3}$ \\ \hline
2& (2,2)&$16$ &$-\frac{8}{3}$ &0&$-264$ \\ \hline
2& (4,1)&$64$ &$-\frac{32}{3}$ &0&$-\frac{2174}{3}$ \\ \hline
2& (6,0)&$320$ &$-\frac{160}{3}$ &0&$-\frac{5956}{3}$ \\ \hline
3&(1,4)& $64$ &$-\frac{56 }{3}$& $0$&$-\frac{104500}{3}$ \\ \hline
3&(3,3)& $320$ &$-\frac{280}{3}$& $0$&$-\frac{429196}{3}$ \\ \hline
3&(5,2)& $2048$ &$-\frac{1792}{3}$& $0$&$-\frac{1759552}{3}$ \\ \hline
3&(7,1)& $15104$ &$-\frac{13216}{3}$& $0$&$-\frac{7209584}{3}$ \\ \hline
3& (9,0)&$123904 $&$-\frac{108416}{3}$&$0$&$-\frac{29527616}{3}$ \\ \hline
4&(0,6)& $384 $ &$-160$& $0$&$-\frac{18667312}{3}$ \\ \hline
4&(2,5)& $2560 $ &$-\frac{3200}{3}$& $0$&$-\frac{101879272}{3}$ \\ \hline
4&(4,4)& $18944$ &$-\frac{22912}{3}$& $256$&$-\frac{555449168}{3}$ \\ \hline
4&(6,3)& $163840$ &$-\frac{194048}{3}$& $3584$&$-\frac{3026251616}{3}$ \\ \hline
4&(8,2)& $1583104$ &$-\frac{1849856}{3}$& $43008$&$-\frac{16485590720}{3}$ \\ \hline
4&(10,1)& $16687104$ &$-6440960 $& $512000$&$-\frac{89806527616}{3}$ \\ \hline
4& (12,0)&$189358080$&$-72652800 $&$6246400$&$-163085218816$ \\ \hline
5&(1,7)& $27136 $ &$-\frac{41792}{3} $& $768$&$-\frac{28726121392}{3}$ \\ \hline
5&(3,6)& $229376$ &$-\frac{331264}{3} $& $13824$&$-\frac{195282001984}{3}$ \\ \hline
5&(5,5)& $2232320$ &$-\frac{3049984}{3} $& $192512$&$-\frac{1326874482304}{3}$ \\ \hline
5&(7,4)& $24391680$ &$-10660352$& $2551808$&$-\frac{9013280450048}{3}$ \\ \hline
5&(9,3)& $291545088$ &$-123583488 $& $34336768$&$-\frac{61226330115584}{3}$ \\ \hline
5&(11,2)& $3750199296$ &$-1553444864$& $477913088$&$-138652119786496$ \\ \hline
5&(13,1)& $51384877056 $ &$-20917362688 $&$6916112384$&$-\frac{2826429058966016}{3}$ \\ \hline
5&(15,0)& $744875950080 $ &$-299359264768 $&$104115208192$&$-\frac{19209989184830464}{3}$ \\ \hline
\end{tabular}
\end{table}

\appendix

\section{Appendix: Proof of Proposition~\ref{Proposition1}}\label{appendixA}

In this appendix, we prove Proposition~\ref{Proposition1}. The proof fundamentally follows the same approach as the corresponding proposition in~\cite{MK}; specifically, we only need to show that the residue of~$f_\Gamma$ for type~(iii) graphs at $w=z_{0}$ vanishes. If $a_i=1$, then $q(x,y)=\prod_{j=1}^{a_i-1}(jx+(a_i-j)y):=1$.

We assume
\[
P(a_1,a_2,\dots,a_N\mid k_1,k_2,\dots,k_m)=P(1,1,\dots,1,a_{r+1},a_{r+2},\dots,a_N\mid k_1,k_2,\dots,k_m),
\]
 which implies $a_{1}=\cdots=a_{r}=1$. We also assume that $a_i > 1$ for $i \in \{r+1,r+2,\dots,N\}$. From the Calabi--Yau condition, we have $r+\sum_{i=r+1}^{N} a_i=\sum_{j=1}^m k_m$.
In the Calabi--Yau case, the condition $n_{a}=0$, $a=2,3,\dots,N-m-1$, holds. Hence, the corresponding integrand is given by\looseness=-1
\begin{gather*}
f_\Gamma =\frac{\text{Sym}(\sigma)}{24\bigl(\prod_{i=1}^Na_i\bigr)^{d}(z_0)^{N(f-1)}} \Biggl( \prod_{i=1}^l\Biggl(\prod_{j=1}^{d_i}\frac{1}{(z_{i,j})^N}
\Biggr)\Biggr)\left(-\frac{N-m}{N}\frac{1}{w^N}-\frac{N+m}{N}\frac{1}{(z_0)^N} \right) \\
\hphantom{f_\Gamma =}{}\times \frac{1}{(w-z_0)^2q(w,z_0)(q(z_0,z_0))^{f-1}} \Biggl( \prod_{p=1}^m\frac{1}{(k_pz_0)^{l-1}}\frac{e_{k_p}(w,z_0) (e_{k_p}(z_0,z_0) )^{f-1}}{k_p w (k_p z_0 )^{f}}\Biggl) \\
\hphantom{f_\Gamma =}{}
\times \Biggl( \prod_{i=1}^l\frac{\prod_{p=1}^me_{k_p}(z_0,z_{i,1})}{q(z_0,z_{i,1})(z_{i,1}-z_0)}\Biggr)\\
\hphantom{f_\Gamma =}{}\times
 \Biggl( \prod_{i=1}^l\Biggl(\prod_{j=1}^{d_i-1}\frac{\prod_{p=1}^m e_{k_p}(z_{i,j},z_{i,j+1})}{q(z_{i,j},z_{i,j+1})(2z_{i,j}-z_{i,j-1}-z_{i,j+1})\bigl(\prod_{p=1}^m k_pz_{i,j} \bigr)} \Biggr) \Biggr) \notag \\
\hphantom{f_\Gamma }{}
=\left(-\frac{N-m}{N}\frac{1}{w^N}-\frac{N+m}{N}\frac{1}{(z_0)^N} \right)\frac{1}{(w-z_0)^2q(w,z_0)}\Biggl( \prod_{l=1}^m\frac{e_{k_l}(w,z_0)}{k_l w}\Biggr)g(z).
\end{gather*}
Here, $g(z)$ in the last line represents the factor that does not contain $w$.

We note the following equalities:
\begin{gather*}
e_k(z_0,z_0)=(kz_0)^{k+1},\qquad
\frac{\mathrm{d}}{\mathrm{d}w}e_k(w,z_0)\Big|_{w=z_0}=\frac{k(k+1)}{2}(kz_0)^k,\\
q(z_0,z_0)=\prod_{i=r+1}^{N}(a_iz_0)^{a_i-1},\\
\frac{\mathrm{d}}{\mathrm{d} w}q(w,z_0)\Big|_{w=z_0}= \sum_{i=r+1}^N \Biggl(\prod_{\substack{l=r+1 \\ l\neq i}}^{N} (a_lz_0)^{a_l-1}\Biggr) (a_iz_0)^{a_i-2}\frac{a_i(a_i-1)}{2}.
\end{gather*}
Using these equalities, the residue of $f_\Gamma$ at $w=z_0$ is computed as follows:
\begin{gather*}
\lim_{w\to z_0} \frac{\mathrm{d}}{\mathrm{d}w}(f_\Gamma)
 =\frac{\mathrm{d}}{\mathrm{d} w}\Biggr[\left(-\frac{N-m}{N}\frac{1}{w^N}-\frac{N+m}{N}\frac{1}{(z_0)^N} \right)\frac{1}{q(w,z_0)}\Biggl( \prod_{l=1}^m\frac{e_{k_l}(w,z_0)}{k_l w}\Biggr)\Biggr]_{w=z_0}g(z) \\
\qquad =g(z)\Biggl[\left(\frac{N-m}{(z_0)^{N+1}} \right)\frac{1}{q(z_0,z_0)}\Biggl( \prod_{l=1}^m\frac{e_{k_l}(z_0,z_0)}{k_l z_0}\Biggr)\\
\qquad\quad{}+\frac{2}{(z_0)^N}\left[\frac{\frac{\mathrm{d}}{\mathrm{d} w}q(w,z_0)}{(q(w,z_0))^2}\right]_{w=z_0}\Biggl( \prod_{l=1}^m\frac{e_{k_l}(z_0,z_0)}{k_l z_0}\Biggr) \\
\qquad\quad{}-\frac{2}{(z_0)^Nq(z_0,z_0)}\sum_{i=1}^m \Biggl( \prod_{\substack{l=1 \\ l\neq i}}^m\frac{e_{k_l}(z_0,z_0)}{k_l z_0}\Biggr)\left(\frac{k_i(k_i+1)}{2k_i z_0}(k_iz_0)^{k_i} -\frac{(k_i z_0)^{k_i+1}}{k_i (z_0)^2}\right)\Biggr]
\\
\qquad =g(z)\Biggl[\frac{N-m}{(z_0)^{N+1}} \frac{1}{\prod_{i=r+1}^{N}(a_iz_0)^{a_i-1}}\Biggl( \prod_{l=1}^m(k_l z_0)^{k_l}\Biggr) \\
\qquad\quad{} +\frac{2}{(z_0)^N\bigl(\prod_{i=r+1}^{N}(a_iz_0)^{a_i-1}\bigr)^2}\Biggl( \prod_{l=1}^m(k_l z_0)^{k_l}\Biggr)\\
\qquad \quad{}\times
\sum_{i=r+1}^{N} \Biggl(\prod_{\substack{l=r+1 \\ l\neq i}}^{N} (a_lz_0)^{a_l-1}\Biggr) (a_iz_0)^{a_i-2}\frac{a_i(a_i-1)}{2} \\
\qquad \quad{} -\frac{2\prod_{l=1}^m(k_l z_0)^{k_l}}{(z_0)^N\prod_{i=r+1}^{N}(a_iz_0)^{a_i-1}}\sum_{i=1}^m \frac{k_i-1}{2 z_0}\Biggr]
\\
\qquad =g(z)\frac{\prod_{l=1}^m(k_l z_0)^{k_l}}{(z_0)^{N}\prod_{i=r+1}^{N}(a_iz_0)^{a_i-1}}\Biggl[\frac{N-m}{z_0}+\sum_{i=r+1}^N \frac{a_i-1}{z_0}-\sum_{i=1}^m \frac{k_i-1}{ z_0}\Biggr]\notag \\
\qquad =g(z)\frac{\prod_{l=1}^m(k_l z_0)^{k_l}}{(z_0)^{N+1}\prod_{i=r+1}^{N}(a_iz_0)^{a_i-1}}\Biggl[N-m+\sum_{i=r+1}^N a_i-(N-r)-\sum_{i=1}^m k_i+m\Biggr]=0.
\end{gather*}

\subsection*{Acknowledgements}

This research is partially supported by JSPS Grant No.~22K03289. The research of M.J.~is also partially supported by JSPS Grant No.~24H00182.

\pdfbookmark[1]{References}{ref}
\LastPageEnding

\end{document}